\documentclass[11pt,letter]{article}

\usepackage[utf8]{inputenc}

\usepackage[T1]{fontenc}
\usepackage{amsmath,amsthm,amsfonts,amssymb}
\usepackage{mathtools}
\usepackage{graphicx}
\usepackage[colorlinks=true, allcolors=blue]{hyperref}
\usepackage{caption}
\usepackage{subcaption}
\usepackage[top=0.9in,bottom=1.1in,left=1.05in,right=1.05in]{geometry}
\usepackage{cleveref}
\usepackage[normalem]{ulem}
\usepackage{xcolor}
\usepackage{booktabs}
\usepackage{cite}

\usepackage[algo2e, noend, noline, ruled, linesnumbered]{algorithm2e}
\DeclareCaptionStyle{ruled}{labelfont=normalfont,labelsep=colon,strut=off} %
\usepackage{algorithm}
\usepackage{algorithmic}

\numberwithin{equation}{section}

\theoremstyle{plain}

\theoremstyle{definition}

\theoremstyle{remark}

\newcommand{\R}{\mathbb{R}}
\newcommand{\E}{\mathbb{E}}
\newcommand{\F}{\mathbb{F}}
\renewcommand{\P}{\mathbb{P}}
\newcommand{\norm}[1]{\left\lVert#1\right\rVert}
\DeclareMathOperator{\tr}{tr}
\newcommand{\dd}{\ensuremath{\mathrm{d}}}
\newcommand{\sigF}{\mathcal{F}}

\DeclareMathOperator*{\argmax}{arg\,max}
\DeclareMathOperator{\unif}{unif}

\newcommand{\revision}[1]{\leavevmode{\color{blue} #1}}

\title{Deep Reinforcement Learning for Infinite Horizon Mean Field Problems in Continuous Spaces}
\author{Andrea Angiuli\thanks{Prime Machine Learning Team, Amazon. 320 West lake Ave N, SEA83, Seattle, WA, 98109, {\em aangiuli@amazon.com}. The work presented here does not relate to this author’s position at Amazon.} \and Jean-Pierre Fouque\thanks{Department of Statistics and Applied Probability, University of California, Santa Barbara, CA 93106-3110, {\em fouque@pstat.ucsb.edu}.} \and Ruimeng Hu\thanks{Department of Mathematics, and Department of Statistics and Applied Probability, University of California, Santa Barbara, CA 93106-3080, {\em rhu@ucsb.edu}.} \and Alan Raydan\thanks{Department of Mathematics, University of California, Santa Barbara, CA 93106-3080, {\em alanraydan@ucsb.edu}. }}

\date{\today}

\begin{document}

\maketitle

\begin{abstract}
We present the development and analysis of a reinforcement learning (RL) algorithm designed to solve continuous-space mean field game (MFG) and mean field control (MFC) problems in a unified manner. The proposed approach pairs the actor-critic (AC) paradigm with a representation of the mean field distribution via a parameterized score function, which can be efficiently updated in an online fashion, and uses Langevin dynamics to obtain samples from the resulting distribution. The AC agent and the score function are updated iteratively to converge, either to the MFG equilibrium or the MFC optimum for a given mean field problem, depending on the choice of learning rates. A straightforward modification of the algorithm allows us to solve mixed mean field control games (MFCGs). The performance of our algorithm is evaluated using linear-quadratic benchmarks in the asymptotic infinite horizon framework.
\end{abstract}


\section{Introduction}
\emph{Mean field games} (MFG) and \emph{mean field control} (MFC)---collectively dubbed mean field problems---are mathematical frameworks used to model and analyze the behavior and optimization of large-scale, interacting agents in settings with varying degrees of cooperation. Since the early 2000s, with the seminal works  \cite{lions2007,malhame2006}, MFGs have been used to study the equilibrium strategies of competitive agents in a large population, accounting for the aggregate behavior of the other agents. Alternately, MFC, which is equivalent to optimal control of McKean-Vlasov SDEs \cite{mckean1966,mckean1967}, focuses on optimizing the behavior of a central decision-maker controlling the population in a cooperative fashion. Cast in the language of stochastic optimal control, both frameworks center on finding an optimal control $\alpha_t$ which minimizes a cost functional objective $J(\alpha)$ subject to given state dynamics in the form of a stochastic differential equation. What distinguishes mean field problems from classical optimal control is the presence of the mean field distribution $\mu_t$, which may influence both the cost functional and the state dynamics. The mean field is characterized by a flow of probability measures that emulates the effect of a large number of participants whose individual states are negligible but whose influence appears in the aggregate. In this setting, the state process $X_t$ models a representative player from the crowd in the sense that the mean field should ultimately be the law of the state process: $\mu_t = \mathcal{L}(X_t)$. The distinction between MFG and MFC, a competitive game versus a cooperative governance, is made rigorous by precisely how we enforce the relationship between $\mu_t$ and $X_t$. We will address the details of the MFG/MFC dichotomy in greater depth in \Cref{sec: mf problems}.

MFG and MFC theories have been instrumental in understanding and solving problems in a wide range of disciplines, such as economics, social sciences, biology, and engineering. In finance, mean field problems have been applied to model and analyze the behavior of investors and markets. For instance, MFG can be used to model the trading strategies of individual investors in a financial market, taking into account the impact of the overall market dynamics. Similarly, MFC can help optimize the management of large portfolios, where the central decision-maker seeks to maximize returns while considering the average behavior of other investors. For in-depth examples of mean field problems in finance, we refer the reader to \cite{carmona2015,carmona-mfg,carmona-lecture}.


Although traditional numerical methods for solving MFG and MFC problems have proceeded along two avenues, solving a pair of coupled partial differential equations (PDE) \cite{carmona2021b} or a forward-backward system of stochastic differential equations (FBSDE) \cite{angiuli2019}, there has been growing interest in solving mean field problems in a model-free way \cite{unified_q_learning,capponi_lehalle_2023,guo2019,perrin2021,lauriere2022,carmona2021}. With this in mind, we turn to \emph{reinforcement learning} (RL), an area of machine learning that trains an agent to make optimal decisions through interactions with a ``black box'' environment. RL can be employed to solve complex problems, such as those found in finance, traffic control, and energy management, in a model-free manner. A key feature of RL is its ability to learn from trial-and-error experiences, refining decision-making policies to maximize cumulative rewards. \emph{Temporal difference} (TD) methods \cite{sutton1988} are a class of RL algorithms that are particularly well-suited for this purpose. They estimate value functions by updating estimates based on differences between successive time steps, combining the benefits of both dynamic programming and Monte Carlo approaches for efficient learning without requiring a complete model of the environment. For a comprehensive overview of the foundations and numerous families of RL strategies, consult \cite{Sutton1998}. \emph{Actor-critic} (AC) algorithms---the modern incarnations of which were introduced in \cite{degris2012}---are a popular subclass of TD methods where separate components, the actor and the critic, are used to update estimates of both a policy and a value function. The actor is responsible for selecting actions based on the current policy, while the critic evaluates the chosen actions and provides feedback to update the policy. By combining the strengths of both policy- and value-based approaches, AC algorithms achieve more stable and efficient learning.

The mean field term itself poses an interesting problem regarding how to numerically store and update a probability measure on a continuous space in an efficient manner. Some authors have chosen to discretize the continuous space, which leads to a vectorized representation as in \cite{unified_q_learning,capponi_lehalle_2023}, while others have looked towards deep learning and deep generative models \cite{perrin2021}. We extend the latter avenue by considering a method of distributional learning known as \emph{score-matching} \cite{hyvarinen2005} in which a probability distribution is represented by the gradient of its log density, i.e., its Stein score function. A parametric representation of the score function is updated using samples from the underlying distribution and allows us to compute new samples from the distribution using a discrete version of Langevin dynamics. We explain how to modify the score-matching procedure for our online regime in \Cref{subsec: score-matching}.

Building off of the work of \cite{unified_q_learning,capponi_lehalle_2023}, in which the authors adapt tabular Q-learning \cite{watkins1989} to solve discrete-space MFG and MFC problems, the main contributions of this papers are:
\begin{enumerate}
\item We address continuous-state and -action space mean field problems in a unified manner by developing an AC algorithm inspired by advantage actor-critic \cite{pmlr-v48-mniha16}. This means that, for a given mean field problem, we use the \emph{same} algorithm to converge to the MFG and MFC solutions simply by adjusting the relative learning rates of the parametric representations of the actor, critic, and mean field distribution. Compared to Q-learning, our proposed algorithm excels in handling continuous or large finite state- and action-spaces, thanks to our neural network implementation. Additionally, it has the advantage of directly learning a policy.

\item  To effectively represent the mean field distribution in continuous space, we employ the \emph{score function}, defined as the gradient of its log density. We parameterize the score function using neural networks and update it via score-matching techniques \cite{hyvarinen2005}. Our method then combines the mean field with the actor-critic paradigm by concurrently learning the score function of the mean field distribution alongside the optimal control, which we derive from the policy learned by the actor.
    
\end{enumerate}

Note that our approach relies on a randomized policy for the actor while another direction of research utilizes entropy-regularization of the value function as in \cite{guo2022, cui2021, JMLR:v21:19-144}. See also the comment in \Cref{subsec: mf hyperperams} describing the relationship between the two approaches. Other papers have concentrated on reinforcement learning for McKean-Vlasov control problems in continuous-time such as in \cite{wei2023}, but our approach treats both MFG and MFC in the same algorithm. We also refer to \cite{frikha2023} for learning MFC problems when the policy may depend additionally on the population distribution.

The rest of the paper is organized as follows. In \Cref{sec: mf problems,sec: rl_ac}, we review the infinite horizon formulation for asymptotic mean field problems and recall the relevant background from RL, respectively. In \Cref{sec: mf_rl} we modify the Markov decision process setting of RL to apply to mean field problems and present our central algorithm. Numerical results and comparisons with benchmark solutions are presented in \Cref{sec: numerical}. As a concluding application, we alter the algorithm in \Cref{sec: mfcg} to apply to \emph{mean field control games} (MFCG), an extension of mean field problems combining both MFG and MFC to model multiple large homogeneous populations where interactions occur not only within each group, but also between groups.

\section{Asymptotic Infinite Horizon Mean Field Problems} \label{sec: mf problems}
In this section, we introduce the framework of mean field games and mean field control in the continuous-time infinite horizon setting. We further emphasize the mathematical distinction between the two classes of mean field problems, highlighting that they yield distinct solutions despite the apparent similarities in their formulation.

In both cases, the mathematical setting is a filtered probability space $(\Omega, \sigF, \F = (\sigF_t)_{t \geq 0}, \P)$ satisfying the usual conditions which supports an $m$-dimensional Brownian motion $(W_t)_{t \geq 0}$. The measurable function $f: \R^d \times \mathcal{P}_2(\R^{d}) \times \R^k \to \R$ is known as the running cost, where $\mathcal{P}(\R^d)$ denotes space of probability measures on $\R^d$ and $\mathcal{P}_s(\R^d)$ is the subspace of $\mathcal{P}(\R^d)$ with finite $s^{th}$-moments.
The constant $\beta > 0$ is a discount factor, which measures the relative importance or value of future rewards or costs compared to immediate ones. For the state dynamics we have drift $b: \R^d \times \mathcal{P}_2(\R^{d}) \times \R^k \to \R^d$ and volatility $\sigma: \R^d \times \mathcal{P}_2(\R^{d}) \times \R^k \to \R^{d \times m}$. 

We will focus on the asymptotic formulation of the infinite horizon mean field problem. In this formulation, we seek a control in the feedback form that depends solely on the state process of the representative player $X: [0, \infty) \times \Omega \to \R^d$, with no explicit time dependency. In other words, the control function is of the form $\alpha: \R^d \to \R^k$, and the trajectory of the control will be given by $\alpha_t = \alpha(X_t)$. This choice allows us to frame the problem more naturally in terms of the Markov decision process setting of reinforcement learning (see \Cref{sec: rl_ac}), which is naturally formulated with time-independent policies.

It was shown in \cite{unified_q_learning} that the asymptotic formulation of the MF problem is the limit of the traditional time-dependent problem as $t \to \infty$ in the following sense: the equilibrium (resp. optimal) control for the time-dependent MFG (resp. MFC) problem converges to the equilibrium (resp. optimal) control for the asymptotic problem as $t \to \infty$.

\subsection{Mean Field Games} \label{subsec: mfg}
The solution of a mean field game, known as a mean field game equilibrium, is a control-mean field pair
\[
    (\hat{\alpha}, \hat{\mu}) \in \mathbb{A} \times  \mathcal{P}(\mathbb{R}^{d}),
\]
where $\mathbb{A}$ is the set of measurable functions $\alpha: \R^d \to \R^k$, satisfying the following conditions:
\begin{enumerate}
    \item $\hat{\alpha}$ solves the stochastic optimal control problem
    \begin{equation} \label{eq: mfg cost}
        \inf _{\alpha \in \mathbb{A}} J_{\hat{\mu}}(\alpha)=\inf _{\alpha \in \mathbb{A}} \mathbb{E}\left[ \int_0^\infty e^{-\beta t} f\left(X_t^{\alpha, \hat{\mu}}, \hat{\mu}, \alpha(X_t^{\alpha, \hat{\mu}})\right)\, \dd t \right],\quad \beta>0,
    \end{equation}
    subject to
    \begin{equation} \label{eq: mfg dynamics}
        \dd X_t^{\alpha, \hat{\mu}}=b\left(X_t^{\alpha, \hat{\mu}}, \hat{\mu}, \alpha(X_t^{\alpha, \hat{\mu}})\right) \, \dd t+\sigma\left(X_t^{\alpha, \hat{\mu}}, \hat{\mu}, \alpha(X_t^{\alpha, \hat{\mu}})\right)\, \dd W_t, \quad X_0^{\alpha, \hat{\mu}}=\xi;
    \end{equation}

    \item $\hat{\mu} = \lim_{t \to \infty} \mathcal{L}(X_t^{\hat{\alpha}, \hat{\mu}})$,
\end{enumerate}
where $\mathcal{L}(X^{\hat{\alpha}, \hat{\mu}}_t)$ refers to the law of $X^{\hat{\alpha}, \hat{\mu}}_t$.

This problem models a scenario in which an infinitesimal player seeks to integrate into a crowd of players already in the asymptotic regime as time tends toward infinity. The resulting stationary distribution represents a Nash equilibrium under the premise that any new player entering the crowd sees no benefit in diverging from this established asymptotic behavior.

\subsection{Mean Field Control} \label{subsec: mfc}
A mean field control problem solution is a control $\alpha^* \in \mathbb{A}$ which satisfies an optimal control problem with McKean-Vlasov dynamics:
\begin{equation} \label{eq: mfc cost}
    \inf _{\alpha \in \mathbb{A}} J(\alpha)=\inf _{\alpha \in \mathbb{A}} \mathbb{E}\left[ \int_0^\infty e^{-\beta t} f\left(X_t^{\alpha}, \mu^\alpha, \alpha(X_t^{\alpha})\right)\, \dd t \right],
\end{equation}
subject to
\begin{equation} \label{eq: mfc dynamics}
     \dd X_t^{\alpha}=b\left(X_t^{\alpha}, \mu^\alpha, \alpha(X_t^{\alpha})\right) \, \dd t+\sigma\left(X_t^{\alpha}, \mu^\alpha, \alpha(X_t^{\alpha})\right)\, \dd W_t, \quad X_0^{\alpha}=\xi,
\end{equation}
using the notation $\mu^\alpha = \lim_{t \to \infty} \mathcal{L}(X_t^{\alpha})$. We will also adopt the notation $\mu^*$ to refer to $\mu^{\alpha^*}$---the limiting distribution for the mean field distribution under the optimal control. In this alternate scenario, we are considering the perspective of a central organizer. Their objective is to identify the control which yields the best possible stationary distribution, ensuring that the societal costs incurred are the lowest possible when a new individual integrates into the group.

Although the initial distribution $\xi$ is specified in both cases, under suitable ergodicity assumptions, the optimal controls $\hat{\alpha}$ and $\alpha^*$ are independent of this initial distribution. For an in-depth treatment of infinite horizon mean field problems, with explicit solutions for the case of linear dynamics and quadratic cost, refer to \cite{malhame2020}.

\subsection{Mean Field Game/Control Distinction}
We summarize the crucial mathematical distinction between MFG and MFC. In the former, one must solve an optimal control problem depending on an arbitrary distribution $\mu$ and then recover the mean field $\hat{\mu}$, which yields the law of the optimal limiting state trajectory. If we consider the map
\[
    \Phi(\mu) = \lim_{t \to \infty} \mathcal{L}(X^{\tilde{\alpha}, \mu}_t),
\]
where $\tilde{\alpha} = \arg \min J_{\mu}(\alpha)$,
then the MFG equilibrium arises as a fixed point of $\Phi$ in the sense that
\[
    \hat{\mu} = \Phi(\hat{\mu}).
\]
In the latter case, the mean field is explicitly the law of the state process throughout the optimization and should be thought of as a pure control problem in which the law of the state process influences the state dynamics. 
Note that in the MFC case, the distribution $\mu^\alpha$ ``moves'' with the choice of control $\alpha$, while in the MFG case, it is ``frozen'' during the optimization step and then a fixed point problem is solved. These interpretations play a key role in guiding the development of this paper's central algorithm, detailed in \Cref{sec: mf_rl}.

Crucially, we conclude this section by noting that, in general,
\begin{equation}
    (\hat{\alpha}, \hat{\mu}) \neq (\alpha^*, \mu^*),
\end{equation}
for the same choice of running cost, discount factor, and state dynamics. Indeed, we will encounter examples of mean field problems with differing solutions when we test our algorithm against benchmark problems in \Cref{sec: numerical}.

\section{Reinforcement Learning and Actor-Critic Algorithms} \label{sec: rl_ac}
Reinforcement learning is a family of machine learning strategies aimed at choosing the sequence of actions which maximizes the long-term aggregate reward from an environment in a model-free way, i.e., assuming no explicit knowledge of the state dynamics or the reward function. Intuitively, one should imagine a black box environment in which an autonomous agent makes decisions in discrete time and receives immediate feedback in the form of a scalar reward signal.

At stage $n$, the agent is in a state $X_{t_n}$ from a given set of states $\mathcal{X}$ and selects an action $A_{t_n}$ from a set of actions $\mathcal{A}$. The environment responds by placing the agent in a new state $X_{t_{n+1}}$ and bestowing it with an immediate reward $r_{t_{n+1}} \in \R$. The agent continues choosing actions, encountering new states, and obtaining rewards in an attempt to maximize the total expected discounted return
\begin{equation}\label{eq: cummul reward}
    \E \left[ \sum_{n=0}^\infty \gamma^n r_{t_{n+1}} \right],
\end{equation}
where $\gamma \in (0,1)$ is a discount factor specifying the degree to which the agent prioritizes immediate reward over long-term returns. The case in which we seek to minimize cost instead of maximize reward as in most financial applications, can be recast in the above setting by taking $r_{t_{n+1}} = -c_{t_{n+1}}$ where $c_{t_{n+1}}$ is the immediate cost incurred at the $n^{th}$ time-step. The expectation in \cref{eq: cummul reward} refers to the stochastic transition from $X_{t_n}$ to $X_{t_{n+1}}$ and, eventually, to the randomness in the choice of $A_{t_n}$. When the new state $X_{t_{n+1}}$ and immediate reward $r_{t_{n+1}}$ only depend on the preceding state $X_{t_n}$ and action $A_{t_n}$, the above formulation is known as a Markov decision process (MDP).

The agent chooses its actions according to a policy $\pi: \mathcal{X} \to \mathcal{P}(\mathcal{A})$, which defines the probability that a certain action should be taken in a given state. The goal of the agent is then to find an optimal policy $\pi^*$ satisfying
\[
    \pi^* \in \argmax_{\pi} \E_{\pi} \left[ \sum_{n=0}^\infty \gamma^n r_{t_{n+1}} \right].
\]
As the reward $r_{t_{n+1}} = r(X_{t_n}, A_{t_n})$ is a function of the current state and current action, the value to be maximized does indeed depend on the policy $\pi$.

For a given policy $\pi$, two quantities of interest in RL are the so-called \emph{state-value function} $v_\pi:\mathcal{X} \to \R$ and the \emph{action-value function} $q_\pi: \mathcal{X} \times \mathcal{A} \to \R$ given by
\begin{align}
    v_\pi(x) &= \E_{\pi} \left[ \sum_{n=0}^\infty \gamma^n r(X_{t_n}, A_{t_n}) \mid X_{t_0} = x \right],\\
    q_\pi(x, a) &= \E_{\pi} \left[ \sum_{n=0}^\infty \gamma^n r(X_{t_n}, A_{t_n}) \mid X_{t_0} = x, A_{t_0} = a \right].
\end{align}
The state-value function defines the expected return obtained from beginning in an initial state $x$ and following the policy $\pi$ from the get-go, whereas the action-value function defines the expected return starting from $x$, taking an initial action $a$, and then proceeding according to $\pi$ after the first step. Moreover, $v_\pi$ and $q_\pi$ are related to each other via the following:
\[
v_{\pi}(x) = \sum_{a \in \mathcal{A}} \pi(a \mid x)  q_\pi(x, a).
\]

The action-value function is integral to many RL algorithms since, assuming that the action-value function $q_*$ corresponding to an optimal policy is known, one can derive an optimal policy by taking the uniform distribution over the actions that maximize $q_*$:
\[
    \pi^*(\cdot \mid  x) = \unif\left(\argmax_{a\in \mathcal{A}} q_*(x, a)\right).
\]
However, since this paper makes far more use of $v_\pi$ than $q_\pi$, we will henceforth refer to the former simply as the ``value function" and the latter as the ``Q-function" as is common in the literature. When referring to both $v_\pi$ and $q_\pi$, we may refer to them jointly as the value functions associated with $\pi$.

For our continuous-space setting, we randomize the policy according to a Gaussian distribution as explained in \Cref{subsec: ac}.

\subsection{Temporal Difference Methods}\label{subsec: td}
In the search for an optimal policy, one often begins with an arbitrary policy, which is improved as the RL agent gains experience in the environment. A key factor in improving a policy is an accurate estimate of the associated value function since this allows us to quantify precisely how much better one policy is over another. 
The value function satisfies the celebrated \emph{Bellman equation}, which relates the value of the current state to that of the successor state:
\begin{equation}\label{eq: bellman}
    v_\pi(x) = \E_\pi [r_{t_{n+1}} + \gamma v_\pi(X_{t_{n+1}}) \mid X_{t_n} = x].
\end{equation}
Since the transition from $X_{t_n}$ to $X_{t_{n+1}}$ is Markovian, \cref{eq: bellman} holds for all $n \geq 0$, not just the initial state. Importantly, the Bellman equation uniquely defines the value function for a given $\pi$, a fact which underlies all algorithms under the umbrella of \emph{dynamic programming}. Solving the Bellman equation for $v_\pi$ is impossible without knowing the reward function and state transition dynamics, so an alternative strategy is needed for our model-free scenario. 

Temporal difference methods center around iteratively updating an approximation $V$ to $v_{\pi}$ in order to drive to zero the TD error $\delta_n$,
\begin{equation} \label{eq: td error}
    \delta_n \coloneqq r_{t_{n+1}} + \gamma V(X_{t_{n+1}}) - V(X_{t_n})
\end{equation}
at each timestep. TD methods use estimates of the value at future times to update the value at the current time, a strategy known as ``bootstrapping''. More importantly, the TD methods we reference here require only the immediate transition sequence $\{X_{t_n}, r_{t_{n+1}}, X_{t_{n+1}}\}$ and no information regarding the MDP model.

\subsection{Actor-Critic Algorithms}\label{subsec: ac}
Actor-critic algorithms form a subset of TD methods in which explicit representations of the policy (the actor) and the value function (the critic) are stored. Often, the representation of the policy is a parametric family of density functions in which the parameters are the outputs of another parametric family of functions, e.g., linear functions, polynomials, and neural networks. In the implementation discussed in \Cref{sec: mf_rl}, our actor is represented by a feedforward neural network which outputs the mean and standard deviation of a normal distribution. The action $A_{t_n}$ is then sampled according to this density. The benefit of a stochastic policy such as this is that it allows for more exploration of the environment so that the agent does not myopically converge to a suboptimal policy.

Since the value function simply outputs a scalar, it may be represented by any sufficiently rich family of real-valued functions. Let
\[
    \Pi_{\psi} \approx \pi \qquad \text{and} \qquad V_{\theta} \approx v
\]
be the parametric approximations of $\pi$ and $v$, both differentiable in their respective parameters $\psi$ and $\theta$. The goal for AC algorithms is to converge to an optimal policy by iteratively updating the actor to maximize the value function and updating the critic to satisfy the Bellman equation. For the critic, this suggests minimizing the following loss function at the $n^{th}$ step:
\begin{equation}
    L_V^{(n)}(\theta) \coloneqq \left( r_{n+1} + \gamma V_\theta(X_{t_{n+1}}) - V_\theta(X_{t_n}) \right)^2 \eqqcolon \delta_n^2.
\end{equation}
Note that the terms inside the square are precisely the TD error $\delta_n$ from \cref{eq: td error}.

The traditional gradient TD update treats the term $y_{t_{n+1}} = r_{n+1} + \gamma V_{\theta}(X_{t_{n+1}})$---known as the \emph{TD target}---as a constant and only considers the term $-V_{\theta}(X_{t_n})$ as a function of $\theta$. This yields faster convergence from gradient descent as opposed to treating both terms as variables in $\theta$ \cite{vanHasselt2012}.
With this in mind, the gradient  of $L_V$ is then
\begin{equation}
    \nabla_{\theta} L_V^{(n)}(\theta) = -2 \delta_n \nabla_{\theta} V_{\theta}(X_{t_n}),
\end{equation}
meaning that, with some learning rate $\rho_V > 0$, we can update the parameters of the critic iteratively using gradient descent:
\begin{equation}
    \theta' = \theta + 2\rho_V \delta_n \nabla_{\theta} V_{\theta}(X_{t_n}).
\end{equation}

Updating the actor, on the other hand, is not so obvious since updating $\psi$ to maximize $V_\theta$ requires somehow computing $\nabla_\psi V_\theta$. Since the connection between $\Pi_\psi$ and $V_\theta$ is not explicit, it is not clear how to compute this gradient a priori. Thankfully, the desired relation comes in the form of the \emph{policy gradient theorem} \cite{sutton1999}, which relates a parameterized policy and its value function via the following:
\begin{equation}\label{eq: pgt}
    \nabla_\psi v_{\Pi_\psi}(x) \propto \mathbb{E}_{\Pi_\psi}\left[q_{\Pi_\psi}(X_{t_n}, A_{t_n}) \nabla_\psi \log \Pi_\psi(A_{t_n} \mid X_{t_n})\right]
\end{equation}
for any initial state $x \in \mathcal{X}$, where $v_{\Pi_{\psi}}$ and $q_{\Pi_{\psi}}$ are the true value functions associated with the parameterized policy $\Pi_{\psi}$.

As a result of the Bellman equation, we have the identity $q_{\pi}(x, a) = \E_{\pi} [ r_{t_{n+1}} + \gamma v_\pi(X_{t_{n+1}})]$ given that $r_{t_{n+1}}$ was the reward obtained by taking the action $a$ in the state $x$. Moreover, adding an arbitrary "baseline" value $\lambda$ to $q_{\Pi_\psi}(X_{t_n}, A_{t_n})$ does not alter the gradient in \cref{eq: pgt} as long as $\lambda$ does not depend on the action $A_{t_n}$. A common baseline value demonstrated to reduce variance and speed up convergence is $\lambda = -v_{\Pi_\psi}(X_{t_n})$ \cite{vanHasselt2012}. With this in mind, we can replace $q_{\Pi_{\psi}}(X_{t_n}, A_{t_n})$ in \cref{eq: pgt} with the TD error $\delta_n = r_{n+1} + \gamma v_{\Pi_\psi}(X_{t_{n+1}}) - v_{\Pi_\psi}(X_{t_n})$ which allows us to reuse $\delta_n$ from its role in updating the critic. As a whole, this suggests the following loss function for the actor:
\begin{equation}
    L_{\Pi}^{(n)}(\psi) \coloneqq -\delta_n \log \Pi_\psi(A_{t_n} \mid X_{t_n})
\end{equation}
For a learning rate $\rho_\Pi > 0$, the gradient descent step would then be
\begin{equation}
    \psi' = \psi + \rho_\Pi \delta_n \nabla_{\psi} \log \Pi_\psi(A_{t_n} \mid X_{t_n}).
\end{equation}

In practical applications, updating the actor and critic in the above fashion at each step generally yields convergence to an optimal policy and value function, respectively. While convergence has been proven in the case of linearly parameterized actor and critic \cite{konda2003}, convergence in the general case is still an open problem. 

\subsubsection{Relative Learning Rates for Actor and Critic} \label{subsec: ac-lr}
Since the gradient descent learning rates play a crucial role in the development of our fundamental algorithm presented in \Cref{sec: mf_rl}, we briefly comment on the choice of learning rates in AC algorithms.

The AC framework alternates between two key steps: refining the critic to accurately approximate the value function associated with the actor's policy---known as \emph{policy evaluation}---and updating the actor to maximize the value returned by the critic---known as \emph{policy improvement}. As the policy improvement step relies on the policy gradient theorem (\cref{eq: pgt}), a sufficiently precise critic is required for its success. Hence, the learning rates for the actor and critic are traditionally chosen such that
\[
    \rho_{\Pi} < \rho_V.
\]
This constraint prompts the critic to learn at a quicker pace compared to the actor, thereby ensuring that the value function from the policy evaluation phase closely aligns with the policy's true value function.

\section{Unified Mean Field Actor-Critic Algorithm for Infinite Horizon} \label{sec: mf_rl}
In this section, we introduce a novel \emph{infinite horizon mean field actor-critic} (IH-MF-AC) algorithm for solving both MFG and MFC problems in continuous-time and continuous-space. Although there have been significant strides in recasting the MDP framework for continuous-time using the Hamiltonian of the associated continuous-time control problem as an analog of the Q-function \cite{JMLR:v23:21-0947,JMLR:v23:21-1387,jia2023qlearning,JMLR:v21:19-144}, we instead take the classical approach of first discretizing the continuous-time problem and then applying the MDP strategies discussed in \Cref{sec: rl_ac}. As our focus is aimed at identifying the stationary solution of the infinite horizon mean field problems, discretizing time does not meaningfully depart from the original continuous-time problem presented in \Cref{sec: mf problems}; While in our ongoing work \cite{angiuli2024} where we tackle the finite horizon regime, the time-discretization must be treated with more care since the mean field becomes a flow of probability distributions parameterized by time, and the optimal control also becomes time-dependent in this context. In the sequel, we will first recast the mean field setting from \Cref{sec: mf problems} as a discrete MDP parameterized by the distribution $\mu$ and then lay out the general procedure of the algorithm before addressing the continuous-space representation of $\mu$ via score functions in \Cref{subsec: score-matching}. \Cref{subsec: unifying} addresses the justification for alternating between the MFG and MFC solutions using the actor, critic, and mean field learning rates.

To begin, we fix a small step size $\Delta t > 0$ and consider the resulting time discretization $(t_0, t_1, t_2, \dots)$ where $t_{n} = n\Delta t$. We then rewrite the cost objectives in \cref{eq: mfg cost,eq: mfc cost} as the Riemann sum
\begin{equation}
    \E \left[ \sum_{n=0}^\infty e^{-\beta t_n} f(X_{t_n}, \mu, A_{t_n}) \Delta t \right]
\end{equation}
and the state dynamics in \cref{eq: mfg dynamics,eq: mfc dynamics} as the Euler-Maruyama approximation
\begin{equation} \label{eq: discrete dynamics}
    X_{t_{n+1}} = X_{t_n} + b(X_{t_n}, \mu, A_{t_n}) \Delta t + \sigma(X_{t_n}, \mu, A_{t_n}) \Delta W_n, \qquad \Delta W_n \sim \mathcal{N}(0,\Delta t).
\end{equation}
This reformulation is directly in correspondence with the MDP setting presented in \Cref{sec: rl_ac}---albeit, parameterized by $\mu$. Observe that $r_{t_{n+1}} = -f(X_{t_n}, \mu, A_{t_n}) \Delta t$, $\gamma = e^{-\beta \Delta t}$, and the state transition dynamics are given by \cref{eq: discrete dynamics}. Note that the model and approximation approach used to derive the reward and next state dynamics are our own choice for simulation purposes but are not prescriptive. The agent sees the environment as a black box and has no knowledge of the details of the problem.

In the style of the AC method described in \Cref{subsec: ac}, our algorithm maintains and updates a policy $\Pi_\psi$ and a value function $V_{\theta}$ which are meant as stand-ins for the control $\hat{\alpha}$ (resp. $\alpha^*$) of the MFG (resp. MFC) and the cost functional $J$, respectively. Both are taken to be feedforward neural networks. The third component is the mean field distribution $\mu$, which is updated simultaneously with the actor and critic at each timestep to approximate the law of $X_t$. The procedure at the $n^{th}$ step is as follows: the agent is in the state $X_{t_n}$ as a result of the dynamics in \cref{eq: discrete dynamics} where $\mu$ is replaced with the current estimate of the mean field $\mu_{n-1}$. The value of $X_{t_n}$ is then used to update the mean field, yielding a new estimate $\mu_n$ (see section \Cref{subsec: score-matching} for details). Using the actor's policy, the agent samples an action $A_{t_n} \sim \Pi_{\psi_n}(\cdot \mid X_{t_n})$ and executes it in the environment. The agent receives a reward which, unbeknownst to the agent, is given by $r_{t_{n+1}} = -f(X_{t_n}, \mu_n, A_{t_n})\Delta t$ (the value $r_{t_{n+1}}$ is known to the agent but the form of $f$ is not). The environment places the agent in a new state $X_{t_{n+1}}$ according to \cref{eq: discrete dynamics} using the distribution $\mu_n$ and the action $A_{t_n}$. The previous steps encapsulate what is meant by \emph{model free} in the sense that the agent has no explicit knowledge of the function $f$ nor the state transition dynamics given by $b$ and $\sigma$; it relies only on the immediate information $r_{t_{n+1}}$ and $X_{t_{n+1}}$. Subsequently, $\Pi_{\psi_n}$ and $V_{\theta_n}$ are updated according to the update rules from \Cref{subsec: ac}.
To mimic the infinite horizon regime, we iterate this procedure for a large number of steps until we achieve convergence to the limiting distribution $\hat{\mu}$ (resp. $\mu^*$) and the equilibrium (resp. optimal) control $\hat{\alpha}$ (resp. $\alpha^*$). The complete pseudocode is presented in \Cref{algo: ihmfac}.

\begin{algorithm}[!tb]
   \caption{\textbf{IH-MF-AC: Infinite Horizon Mean Field Actor-Critic}}
   \label{algo: ihmfac}
\begin{algorithmic}[1] 
    \REQUIRE Number of time steps $N \gg 0$; discrete time step size $\Delta t$; neural network learning rates for actor $\rho_\Pi$, critic $\rho_V$, and score $\rho_\Sigma$; Langevin dynamics step size $\epsilon$.
    \STATE Initialize neural networks:\\
    \textbf{Actor} $\Pi_{\psi_0}: \R^d \to \mathcal{P}(\R^k)$\\
    \textbf{Critic} $V_{\theta_0}: \R^d \to \R$\\
    \textbf{Score} $\Sigma_{\varphi_0}: \R^d \to \R^d$\vspace{0.2cm}
    \STATE Agent receives initial state $X_{t_0}$ from the Environment.\vspace{0.2cm}
   \FOR{$n=0,\dots,N-1$}\vspace{0.2cm}
      \STATE Environment computes score loss:
      $\quad L_\Sigma^{(n)} (\varphi_n) = \tr\left( \nabla_x \Sigma_{\varphi_n}(X_{t_n}) \right) + \frac{1}{2}\norm{\Sigma_{\varphi_n}(X_{t_n})}_2^2$\vspace{0.2cm}
      
      \STATE Environment updates score with SGD:
      $\quad \varphi_{n+1} = \varphi_n -\rho_\Sigma \nabla_{\varphi} L_\Sigma^{(n)} (\varphi_n)$\vspace{0.2cm}
      
      \STATE Environment generates mean field samples $S_{t_n} = \left(S_{t_n}^{(1)}, S_{t_n}^{(2)}, \dots, S_{t_n}^{(k)}\right)$ from $\Sigma_{\varphi_{n+1}}$ using Langevin dynamics (\cref{eq: langevin}) with step size $\epsilon$ and compute $\overline{\mu}_{S_{t_n}} \coloneqq \frac{1}{k} \sum_{i=1}^k \delta_{S_{t_n}^{(i)}}$.\vspace{0.2cm}

      \STATE Agent samples action:
      $\quad A_{t_n} \sim \Pi_{\psi_n}(\cdot \mid X_{t_n})$\vspace{0.2cm}
      
      \STATE Agent observes reward $r_{n+1}$ and next state $X_{t_{n+1}}$ generated from the environment based on its knowledge of $\overline{\mu}_{S_{t_n}}$. \vspace{0.2cm}

      \STATE Agent computes TD target:
      $\quad y_{n+1} = r_{n+1} + e^{-\beta \Delta t} V_{\theta_n}(X_{t_{n+1}})$\vspace{0.2cm}
      
      \STATE Agent computes TD error:
      $\quad \delta_{\theta_n} =y_{n+1} - V_{\theta_n}(X_{t_n})$\vspace{0.2cm}
      
      \STATE Agent computes critic loss:
      $\quad L_V^{(n)}(\theta_n) = \delta_{\theta_n}^2$\vspace{0.2cm}
      
      \STATE Agent updates critic with SGD:
      $\quad \theta_{n+1} = \theta_n - \rho_V \nabla_{\theta} L_V^{(n)}(\theta_n)$\vspace{0.2cm}
      
      \STATE Agent computes actor loss:
      $\quad L_{\Pi}^{(n)}(\psi_n) = -\delta_{\theta_n} \log \Pi_{\psi_n}(A_{t_n} \mid X_{t_n})$\vspace{0.2cm}
      
      \STATE Agent updates actor with SGD:
      $\quad \psi_{n+1} = \psi_n - \rho_{\Pi} \nabla_{\psi} L_\Pi^{(n)}(\psi_n)$\vspace{0.2cm}
   \ENDFOR \vspace{0.2cm}
   \RETURN $(\Pi_{\psi_N}, \Sigma_{\varphi_N})$
\end{algorithmic}
\end{algorithm}

\subsection{Representation of $\mu$ via Score-matching}\label{subsec: score-matching}
The question of how to represent and update $\mu$ in the continuous-space setting deserves special consideration in this work. In \cite{capponi_lehalle_2023,unified_q_learning}, the authors deal with the discrete-space mean field distribution in a natural way, using a normalized vector containing the probabilities of each state. Each individual state is modeled as a one-hot vector (a Dirac delta measure), and the approximation $\mu_n$ is updated at each step using an exponentially weighted update of the form $\mu_{n+1} = \mu_{n} + \rho_\mu (\delta_{X_{t_{n}}} - \mu_{n})$ with the mean field learning rate $\rho_\mu > 0$. \cite{frikha2023} uses a similar update in the context of an AC algorithm for solving only MFC problems, while focusing on a more in-depth treatment of the continuous-time aspect. The authors in \cite{perrin2021} tackle continuous-state spaces for the MFG problem using the method of normalizing flows, which pushes forward a fixed latent distribution, such as a Gaussian, using a series of parameterized invertible maps \cite{rezende2015}. There is reason to believe that other deep generative models, such as generative adversarial networks (GANs) or variational auto-encoders (VAEs), may yield successful representations of the population distribution with their own drawbacks and advantages.

In our case, partly due to its simplicity of implementation, we opt for the method known as \emph{score-matching} \cite{hyvarinen2005}, which has been successfully applied to generative modeling \cite{song2019}. If $\mu$ has a density function $p_\mu: \R^d \to \R$, then its Stein score function is defined as
\[
    s_\mu(x) \coloneqq \nabla \log p_\mu(x).
\]
The score function is a useful proxy for $\mu$ in the sense that we can use $s_\mu$ to generate samples from $\mu$ using a Langevin Monte Carlo approach. Given an initial sample $x_0$ from an arbitrary distribution and a small step size $\epsilon > 0$, the sequence defined by
\begin{equation}\label{eq: langevin}
    x_{m+1} = x_m + \frac{\epsilon}{2} s_\mu(x_m) + \sqrt{\epsilon}\,z_m, \qquad z_m \sim \mathcal{N}(0,1)
\end{equation}
converges to a sample from $\mu$ as $m \to \infty$.

From the standpoint of parametric approximation, if $(\Sigma_\varphi)_{\varphi \in \Phi}$ is a sufficiently rich family of functions from $\R^d \to \R^d$, the natural goal is to find the parameters $\varphi$ which minimize the residual $\E_{x \sim \mu} [ \norm{\Sigma_{\varphi}(x) - s_\mu(x)}_2^2 ]$. Although we do not know the true score function, a suitable application of integration by parts yields an expression that is proportional to the previous residual but independent of $s_\mu$:
\begin{equation}\label{eq: expected score loss}
    \E_{x \sim \mu} \left[ \tr(\nabla_x \Sigma_{\varphi}(x)) + \frac{1}{2} \norm{\Sigma_\varphi(x)}^2_2 \right].
\end{equation}

We adapt the above expression for our online setting in the following way. At the $n^{th}$ step, we have a sample $X_{t_n}$ of the state process and a score representation $\Sigma_{\varphi_n}$. We take the loss function for $\Sigma$ to be
\begin{equation}
    L_\Sigma^{(n)} (\varphi_n) \coloneqq \tr\left( \nabla_x \Sigma_{\varphi_n}(X_{t_n}) \right) + \frac{1}{2}\norm{\Sigma_{\varphi_n}(X_{t_n})}_2^2.
\end{equation}
Assuming $\Sigma$ is differentiable with respect to $\varphi$, we then update the parameters using the gradient descent step
\begin{equation}\label{eq: score update}
    \varphi_{n+1} = \varphi_n - \rho_\Sigma \nabla_{\varphi} L_\Sigma^{(n)} (\varphi_n)
\end{equation}
where $\rho_\Sigma > 0$ is the mean field learning rate. Now we can generate samples from $\Sigma_{\varphi_{n+1}}$ and take $\mu_n$ to be the empirical distribution of these samples. More concretely, let $S_{t_n} = \left(S_{t_n}^{(1)}, S_{t_n}^{(2)}, \dots, S_{t_n}^{(k)}\right)$ be the $k$ samples generated from $\Sigma_{\varphi_{n+1}}$ using the Langevin Monte Carlo algorithm in \cref{eq: langevin}, and let
\[
    \mu_n = \overline{\mu}_{S_{t_n}},
\]
where the notation $\overline{\mu}_S \coloneqq \frac{1}{k} \sum_{i=1}^k \delta_{S^{(i)}}$ denotes the empirical distribution of the points $S = (S^{(1)}, S^{(2)}, \dots, S^{(k)})$. By the law of large numbers, $\overline{\mu}_{S_{t_n}}$ converges to the true distribution corresponding to $\Sigma_{\varphi_{n+1}}$ as $k \to \infty$.

In the context of generative modeling, the gradient descent update in \cref{eq: score update} is usually evaluated with several mini-batches of independent samples all from a single distribution. This contrasts with our online approach in which each update is done with the current state $X_{t_n}$, which is generated from a different distribution than the previous state. We justify this as a form of bootstrapping in which we attempt to learn a target distribution that is continuously moving, but ultimately converging to the limiting distribution of the MFG or MFC. Since our updates depend on individual samples, we expect the loss $L_\Sigma$ to be a noisy estimate of the expectation in \cref{eq: expected score loss}, which may slow down convergence.
Rather than updating at every timestep, another option would be to perform a batch update after every $m > 1$ timesteps using all samples $(X_{t_n}, X_{t_{n+1}}, \dots, X_{t_{n+(m-1)}})$ generated along the state trajectory, which may accelerate convergence by reducing variance. It is important to acknowledge that the $m$ samples will come from different distributions, so the batch update will also introduce bias into the gradient estimate. This may be mitigated by instead running multiple trajectories in parallel and updating the score function at each step using the samples $(X^{(1)}_{t_n}, X^{(2)}_{t_n}, \dots, X^{(m)}_{t_n})$ from the same timestep.

\subsection{Unifying Mean Field Game and Mean Field Control Problems}\label{subsec: unifying}
Having laid out the general algorithm, we now address the issue of unifying the MFG and MFC formulations in the style of \cite{unified_q_learning,capponi_lehalle_2023}.

The intuitions presented in \Cref{sec: mf problems} regarding the difference between MFG and MFC suggest that the interplay between the learning rates $\rho_{\Pi}$, $\rho_{V}$, and $\rho_{\Sigma}$ may be used to differentiate between the two solutions of the mean field problem. Taking $\rho_\Sigma < \min\{\rho_\Pi, \rho_V\}$ emulates the notion of solving the classical control problem corresponding to a fixed (frozen)  $\mu$---or, in this case, a slowly moving $\mu$---and then updating the distribution to match the law of the state process in an iterative manner. This matches the strategy discussed in \Cref{subsec: mfg} for finding an MFG equilibrium. Conversely, taking $\rho_\Sigma > \max\{\rho_\Pi, \rho_V\}$ is more in keeping with simultaneous optimization of the mean field and the policy, which should yield the MFC solution as discussed in \Cref{subsec: mfc}.

Borkar's stochastic multi-scale approximation framework \cite{borkar1997,borkar2008} has been applied in \cite{unified_q_learning,capponi_lehalle_2023,angiuli2023convergence} to describe how the choice of learning rates is crucial to derive the solution of different mean field problems using a 2-time scale Q-learning algorithm. Following similar ideas, we show how our method can be expressed as a 3-time scale actor-critic algorithm converging to MFG or MFC solution depending on the choice of learning rates.

In particular, \Cref{algo: ihmfac} is characterized by the following system of numerical updates  
\begin{equation}
\label{eq: discrete system of ODEs}
 \begin{cases}
    \quad \varphi_{n+1} & = \varphi_n -\rho_n^\Sigma \nabla_{\varphi} L_\Sigma (\varphi_n), \\
    \quad \theta_{n+1} & = \theta_n - \rho_n^V \nabla_{\theta} L_V(\theta_n), \\
    \quad \psi_{n+1} & = \psi_n - \rho_n^{\Pi} \nabla_{\psi} L_\Pi(\psi_n),
\end{cases}   
\end{equation}
where $L_\Sigma$, $L_V$, and $L_\Pi$ (each of which actually depend implicitly on all three of the parameter vectors) are the loss functions of the Stein score, the critic, and the actor networks, respectively, with $\rho_n^\Sigma$, $\rho_n^V$ and $\rho_n^{\Pi}$ being the corresponding learning rates satisfying usual Robbins–Monro type conditions, namely: $\sum_{n=0}^{\infty} \rho_n^{\Sigma}=\sum_{n=0}^{\infty} \rho_n^{V}=\sum_{n=0}^{\infty} \rho_n^{\Pi}=\infty$ and $\sum_{n=0}^{\infty} {(\rho_n^{\Sigma})}^2=\sum_{n=0}^{\infty} {(\rho_n^{V})}^2=\sum_{n=0}^{\infty} {(\rho_n^{\Pi})}^2<\infty$ . 

The relationship between  learning rates is crucial to determine the convergence of our algorithm. As shown by \cite{borkar1997actor}, the actor-critic paradigm can be represented as a multi-scale stochastic approximation procedure in which the updates of the actor evolve at a much slower pace than the ones of the critic. Indeed, as discussed in \Cref{subsec: ac-lr}, choosing $\rho_{\Pi}<\rho_V$ allows the actor $\Pi$ to be seen as frozen by the critic $V$ which can be learned accordingly. Our algorithm extends the traditional 2-time scale framework of actor-critic algorithms to approach mean field problems by adding an additional time scale for learning the mean field distributions. 

\paragraph{Three time scale approach for MFG.}
If $\rho^{\Sigma}_n<\rho^{\Pi}_n<\rho^{V}_n$ so that $\rho^{\Sigma}_n/\rho^{\Pi}_n \mapsto 0$ and $\rho^{\Pi}_n/\rho^{V}_n \mapsto 0$ as $n\mapsto \infty$, the system \ref{eq: discrete system of ODEs} tracks the ODE system:

\begin{equation}
 \begin{cases}
    \quad \dot \varphi & =  - \nabla_{\varphi} L_\Sigma (\varphi, \psi, \theta) \\
    \quad \dot \psi & =  -  \frac{1}{\epsilon } \nabla_{\psi} L_\Pi(\varphi, \psi, \theta)\\
    \quad \dot \theta & = - \frac{1}{\epsilon \tilde \epsilon} \nabla_{\theta} L_V(\varphi, \psi, \theta) 
\end{cases}   
\end{equation}
where $\rho^{\Sigma}_n / \rho^{\Pi}_n$ and $\rho^{\Pi}_n / \rho^{V}$ are thought of being of order $\epsilon \ll 1$ and $\tilde\epsilon \ll 1$, respectively. In the case $\varphi$ and $\psi$ are fixed, we assume the ODE $\dot \theta  = - \frac{1}{\epsilon \tilde \epsilon} \nabla_{\theta} L_V(\varphi, \psi, \theta)$ to have a stable equilibrium $(\varphi, \psi, \theta^{\varphi,\psi})$ corresponding to the local minimum of the operator $L_V$, meaning that the network approximates the value function $V$ corresponding to the policy $\Pi$ and the population distribution $\mu$ described respectively by the parameters $\psi$ and $\varphi$. Similarly, given fixed $\varphi$, we assume the ODE $\dot \psi  =  -  \frac{1}{\epsilon } \nabla_{\psi} L_\Pi(\varphi, \psi, \theta_{\varphi,\psi})$ to reach a stable equilibrium $(\varphi, \psi_{\varphi}, \theta_{\varphi, \psi_{\varphi}})$ as the local minimum of the operator $L_{\Pi}$, meaning that the network approximates the optimal policy of the infinitesimal player facing the crowd distribution $\mu$ with Stein score function parameterized by $\varphi$. Finally, the ODE $\dot \varphi =  - \nabla_{\varphi} L_\Pi(\varphi, \psi_{\varphi}, \theta_{\varphi, \psi_{\varphi}})$ stabilizes at $(\varphi^*, \psi_{\varphi^*}, \theta_{\varphi^*, \psi_{\varphi^*}})$, the local minimum of $L_{\Sigma}$ which results the network approximating the Stein score function of the population at equilibrium of the MFG problem.

\paragraph{Three time scale approach for MFC.}
If $\rho^{\Pi}_n<\rho^{V}_n<\rho^{\Sigma}_n$ so that $\rho^{\Pi}_n/\rho^{V}_n \mapsto 0$ and $\rho^{V}_n/\rho^{\Sigma}_n \mapsto 0$ as $n\mapsto \infty$, the system \ref{eq: discrete system of ODEs} tracks the ODE system:
\begin{equation}
 \begin{cases} 
    \quad \dot \psi & =  - \nabla_{\psi} L_\Pi(\varphi, \psi, \theta) \\
    \quad \dot \theta & = - \frac{1}{\epsilon } \nabla_{\theta} L_V(\varphi, \psi, \theta) \\
    \quad \dot \varphi & =  - \frac{1}{\epsilon \tilde \epsilon} \nabla_{\varphi} L_\Sigma (\varphi, \psi, \theta)
\end{cases}   
\end{equation}
where $\rho^{\Pi}_n / \rho^{V}_n$ and $\rho^{V}_n / \rho^{\Sigma}$ are thought of being of order $\epsilon \ll 1$ and $\tilde\epsilon \ll 1$ respectively. Considering $\psi$ and $\theta$ to be fixed, we assume the ODE $\dot \varphi =  - \frac{1}{\epsilon \tilde \epsilon} \nabla_{\varphi} L_\Sigma (\varphi, \psi, \theta)$ to have a stable equilibrium $( \varphi_{\psi, \theta}, \psi, \theta)$ corresponding to the local minimum of the operator $L_\Sigma$, meaning that the network approximates the population distribution induced by the policy and value function described respectively by the parameter $\psi$ and $\theta$. Similarly, fixed $\psi$, we assume the ODE $\dot \theta = - \frac{1}{\epsilon } \nabla_{\theta} L_V(\varphi_{\psi, \theta}, \psi, \theta)$ to reach a stable equilibrium $(\varphi_{\psi, \theta_{\psi}}, \psi, \theta_{\psi})$ as the local minimum of the operator $L_V$, meaning that the network approximates the value function corresponding to the policy $\Pi$ described  by $\psi$. Finally, the ODE $\dot \psi = - \nabla_{\psi} L_\Pi(\varphi_{\psi, \theta_{\psi}}, \psi, \theta_{\psi}))$ stabilizes at $(\varphi_{\psi^*, \theta_{\psi^*}}, \psi^*, \theta_{\psi^*})$, the local minimum of $L_{\Pi}$ which corresponds to the network approximating the optimal control strategy of the MFC problem.

Writing a rigorous proof of convergence with precise assumptions on these operators is extremely complicated and outside of the scope of this paper. It was recently achieved in the case of Q-learning with finite-state and -action spaces in \cite{angiuli2023convergence}.

\section{Numerical Results} \label{sec: numerical}

\subsection{A Linear-Quadratic Benchmark}\label{subsec: lq bench}
We test our algorithm on a linear-quadratic (LQ) mean field problem where we wish to optimize
\begin{equation} \label{eq: lq cost}
    \E \left[ \int_0^\infty e^{-\beta t} \left(\frac{1}{2} \alpha_t^2+c_1\left(X_t-c_2 m\right)^2+c_3\left(X_t-c_4\right)^2+c_5 m^2\right) \, \dd t \right]
\end{equation}
with state dynamics
\begin{equation} \label{eq: lq dynamics}
    \dd X_t = \alpha_t\, \dd t + \sigma \, \dd W_t, \qquad t \in [0, \infty)
\end{equation}
where $m = \int x \,\mu(\dd x)$ so that the mean field dependence is only through the first moment of the asymptotic distribution $\mu$.
Note that the state dynamics depend only linearly on the control $\alpha$, and the running cost function depends on $\alpha$, $X$, and $m$ quadratically, hence the name linear-quadratic.

The various terms in \cref{eq: lq cost} have the following interpretations: the first and last terms penalize $\alpha$ and $m$ from being too large, the second term addresses the relationship between the state process and the mean field distribution, which penalizes $X$ from deviating too far from $c_2m$, and the third term penalizes $X$ for being far from $c_4$. The coefficients $c_1$, $c_3$, and $c_5$ determine the relative influence of each term on the total cost.

Both of the 1-dimensional MFG and MFC problems corresponding to \cref{eq: lq cost,eq: lq dynamics} have explicit analytic solutions, which we state now using the notation consistent with the derivations in \cite{unified_q_learning}. We also test our algorithm on the equivalent multivariate problem for which we include the full derivation of the solutions in \Cref{ap: derivations}.

\subsection{Solution for Asymptotic Mean Field Game}
Traditional methods for deriving the LQ problem solution begin with recovering the value function
\[
    v(x) \coloneqq \inf _{\alpha \in \mathbb{A}} \mathbb{E}\left[\left.\int_0^{\infty} e^{-\beta t}\left(\frac{1}{2} \alpha_t^2+c_1\left(X_t-c_2 m\right)^2+c_3\left(X_t-c_4\right)^2+c_5 {m}^2\right)\, \mathrm{d}t \right\rvert\, X_0=x\right]
\]
as the solution of a Hamilton-Jacobi-Bellman equation. In the MFG case, we denote the optimal value function as $\hat{v}$ and an explicit form is given by $\hat{v}(x) = \hat{\Gamma}_2 x^2 + \hat{\Gamma}_1 x + \hat{\Gamma}_0$, where
\begin{align*}
    &\hat{\Gamma}_2 =\frac{-\beta+\sqrt{\beta^2+8\left(c_1+c_3\right)}}{4}\\[.5em]
    &\hat{\Gamma}_1 = - \frac{2 \hat{\Gamma}_2 c_3 c_4}{\hat{\Gamma}_2 ( \beta + 2 \hat{\Gamma}_2 ) - c_1 c_2}\\[.5em]
    &\hat{\Gamma}_0 = \frac{c_5 \hat{m}^2+c_3 c_4^2+c_1 c_2^2 \hat{m}^2+\sigma^2 \hat{\Gamma}_2-\frac{1}{2} \hat{\Gamma}_1^2}{\beta}.
\end{align*}
Then the optimal control for the MFG is
\begin{equation} \label{eq: mfg control}
    \hat{\alpha}(x) = -\hat{v}'(x) = -\left(2 \hat{\Gamma}_2 x + \hat{\Gamma}_1\right).
\end{equation}
Substituting \cref{eq: mfg control} into \cref{eq: lq dynamics} yields the Ornstein-Uhlenbeck process
\[
    \dd \hat{X}_t = -\left( 2\hat{\Gamma}_2 \hat{X}_t + \hat{\Gamma}_1 \right)\, \dd t + \sigma \, \dd W_t,
\]
whose limiting distribution $\hat{\mu} = \lim_{t \to \infty} \mathcal{L}(\hat{X}_t)$ is
\begin{equation}
    \hat{\mu} = \mathcal{N}\left( -\frac{\hat{\Gamma}_1}{2 \hat{\Gamma}_2}, \frac{\sigma^2}{4 \hat{\Gamma}_2} \right).
\end{equation}
Since the mean field interaction for the LQ problem is only through the mean $\hat{m} = \int x \,\hat{\mu}(\dd x)$, we note that a simplified form of $\hat{m}$ is
\begin{equation} \label{eq: mfg mean}
    \hat{m} = -\frac{\hat{\Gamma}_1}{2 \hat{\Gamma}_2} = \frac{c_3 c_4}{c_1 + c_3 - c_1 c_2}.
\end{equation}

\subsection{Solution for Asymptotic Mean Field Control}
Similarly as above, we denote the MFC value function by $v^*$ and claim that is has the form $v^*(x) = \Gamma^*_2 x^2 + \Gamma_1^* x + \Gamma^*_0$ with
\begin{align*}
    &\Gamma^*_2 =\frac{-\beta+\sqrt{\beta^2+8\left(c_1+c_3\right)}}{4}\\[.5em]
    &\Gamma^*_1 = - \frac{2 \Gamma_2^* c_3 c_4}{\Gamma^*_2(\beta + 2 \Gamma^*_2) + c_5 - c_1 c_2(2 - c_2)}\\[.5em]
    &\Gamma^*_0 = \frac{c_5 {m^*}^2+c_3 c_4^2+c_1 c_2^2 {m^*}^2+\sigma^2 \Gamma^*_2-\frac{1}{2} {\Gamma^*_1}^2}{\beta}.
\end{align*}
The optimal control for the MFC is
\begin{equation} \label{eq: mfc control}
    \alpha^*(x) = - {v^*}'(x) =  -\left(2 \Gamma_2^* x + \Gamma^*_1\right).
\end{equation}
Substituting \cref{eq: mfc control} into \cref{eq: lq dynamics} yields the Ornstein-Uhlenbeck process
\[
    \dd X^*_t = -\left(2 \Gamma_2^* X^*_t + \Gamma^*_1\right)\, \dd t + \sigma \, \dd W_t,
\]
whose limiting distribution $\mu^* = \lim_{t \to \infty} \mathcal{L}(X^*_t)$ is
\begin{equation}
    \mu^* = \mathcal{N}\left( -\frac{\Gamma^*_1}{2 \Gamma^*_2}, \frac{\sigma^2}{4 \Gamma^*_2} \right).
\end{equation}
Since the mean field interaction is only through the mean $m^* = \int x \, \mu^*(\dd x)$, we note that an equation for $m^*$ which only depends explicitly on the running cost coefficients is
\begin{equation}\label{eq: mfc mean}
    m^* = -\frac{\Gamma^*_1}{2 \Gamma^*_2} = \frac{c_3 c_4}{c_1 +c_3 + c_5 - c_1 c_2 (2- c_2)}.
\end{equation}

\subsection{A multivariate Linear-Quadratic Benchmark}
We also test our algorithm on a higher-dimensional generalization of the scalar LQ problem in \cref{eq: lq dynamics,eq: lq cost} with the cost functional
\begin{equation}\label{eq: multivar cost}
    \mathbb{E} \left[\int_0^{\infty} e^{\beta t} \left(\frac{1}{2} \alpha_t^T \alpha_t + (x - C_2 m)^T C_1 (x - C_2 m) + (x - c_4)^T C_3 (x - c_4) + m^T C_5 m \right) \, \dd t \right],
\end{equation}
and state dynamics
\begin{equation}\label{eq: multivar dynamics}
    \dd X_t = \alpha_t\, \dd t + \sigma \, \dd W_t, \qquad t \in [0, \infty).
\end{equation}
In this case, $C_1, C_3, C_5 \in \R^{d \times d}$ are positive-definite matrices, $c_4 \in \R^d$, $C_2 \in \R^{d \times d}$, $\sigma \in \R^{d \times m}$, and $W_t$ is an $m$-dimensional Brownian motion. We test our algorithm on a benchmark problem for the case $d = m = 2$, and refer the reader to \Cref{ap: derivations} for the derivations of the analytic solutions to both the MFG and MFC problems associated with the above system.

 \subsection{Hyperparameters and Numerical Specifics} \label{subsec: mf hyperperams}
\begin{figure}[!htb]
    \centering
    \begin{minipage}{\textwidth}

        \begin{table}[H]
        \caption{Running cost coefficients and volatility for \cref{eq: lq cost,eq: lq dynamics}. The results for this parameter set are displayed in \cref{fig: results 1,fig: mean errors 1,fig: values 1,fig: value errors 1,fig: control errors 1}.}
        \centering
            \begin{tabular}{cccccc}
            \toprule
                 $c_1$ & $c_2$ & $c_3$ & $c_4$ & $c_5$ & $\sigma$ \\
            \midrule
                 0.25 & 1.5 & 0.5 & 0.6 & 1.0 & 0.3\\
            \bottomrule
            \end{tabular}
            \label{tab: coeffs 1}
        \end{table}

        \vspace{0.5cm}
        
        \begin{minipage}{\textwidth}
              \centering
              \begin{subfigure}{.5\textwidth}
              \centering
                  \includegraphics[width=.99\linewidth]{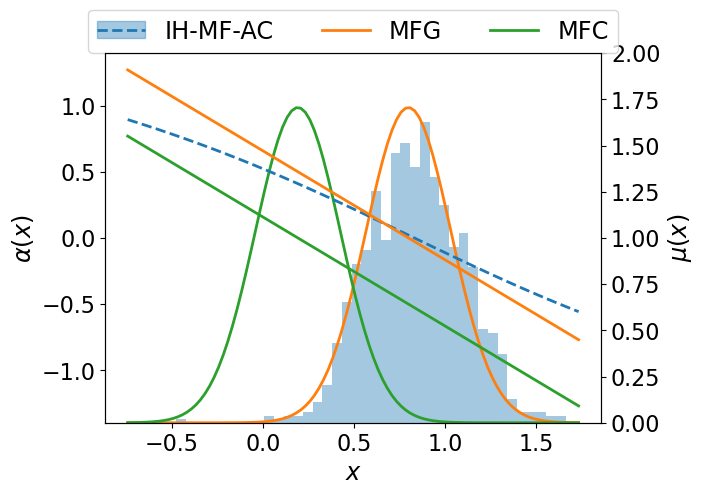}
             \caption{Result for $\rho_\Sigma = 10^{-6}$}
              \label{fig: mfg results1}
              \end{subfigure}%
            \begin{subfigure}{.5\textwidth}
              \centering
              \includegraphics[width=.99\linewidth]{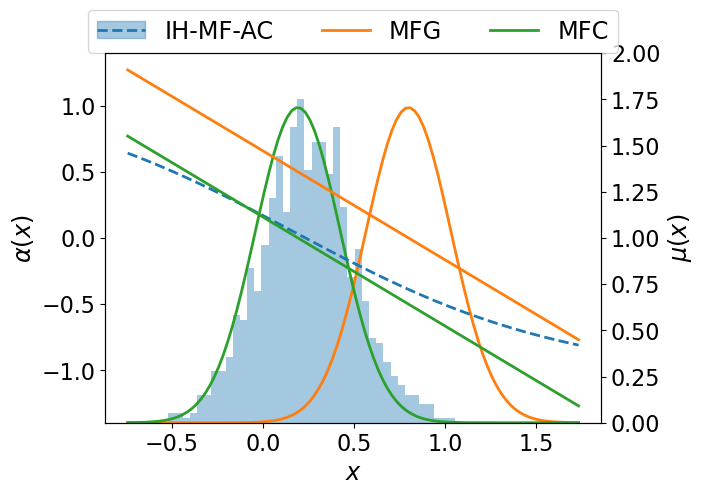}
              \caption{Result for $\rho_\Sigma = 5 \times 10^{-4}$}
              \label{fig: mfc results1}
              \end{subfigure}
              \caption{The histogram (blue) is the learned asymptotic distribution using samples generated from the parameterized score function $\Sigma_{\varphi_N}$ and the dashed line (blue) is the learned feedback control after $N=10^6$ iterations averaged over five runs with different initial samples. The green curves correspond to the optimal control and mean field distribution for MFC, while the orange curves are the equivalent for MFG. The bottom axis shows the state variable $x$, the left axis refers to the value of the control $\alpha(x)$, and the right axis represents the probability density of $\mu$.}
              \label{fig: results 1}
        \end{minipage}

        \vspace{0.5cm}

        \begin{minipage}{\textwidth}
              \centering
              \begin{subfigure}{.5\textwidth}
              \centering
                  \includegraphics[width=.98\linewidth]{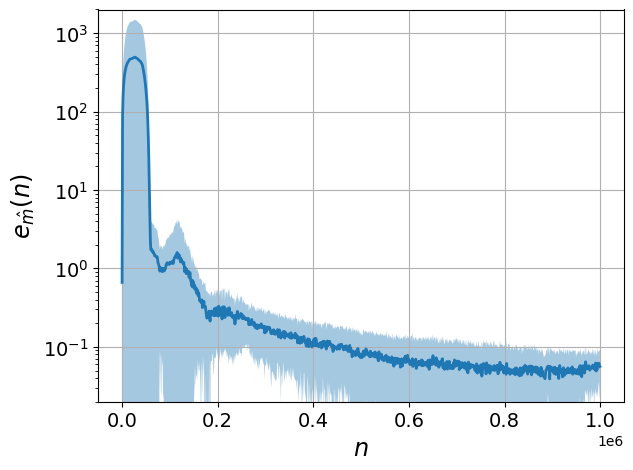}
                  \caption{Results for $\rho_\Sigma = 10^{-6}$}
              \label{fig: mfg mean error 1}
              \end{subfigure}%
            \begin{subfigure}{.5\textwidth}
              \centering
              \includegraphics[width=.98\linewidth]{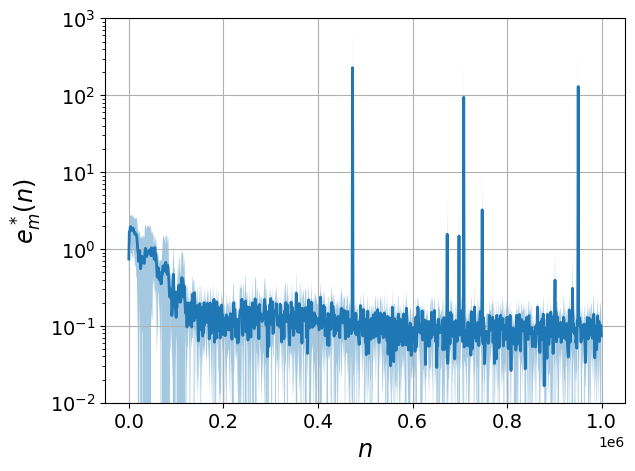}
              \caption{Results for $\rho_\Sigma = 5 \times 10^{-4}$}
              \label{fig: mfc mean error 1}
              \end{subfigure}
               \caption{We plot the absolute error between the mean of samples produced from the parameterized score function $\Sigma_{\varphi_n}$ and the optimal mean $\hat{m}$ in the case of MFG (left) and $m^*$ in the case of MFC (right). These values were averaged over five runs each with different random initial samples with the standard deviation given by the light blue shaded region. Large jumps are due to random outliers which result from the stochasticity of our algorithm.}
              \label{fig: mean errors 1}
        \end{minipage}
    \end{minipage}
\end{figure}

\begin{figure}[!htb]
    \centering
    \begin{minipage}{\textwidth}

        \begin{table}[H]
        \caption{Running cost coefficients and volatility for \cref{eq: lq cost,eq: lq dynamics}. The results for this parameter set are displayed in \cref{fig: results 2,fig: mean errors 2,fig: values 2,fig: value errors 2,fig: control errors 2}.}
        \centering
            \begin{tabular}{cccccc}
            \toprule
                 $c_1$ & $c_2$ & $c_3$ & $c_4$ & $c_5$ & $\sigma$ \\
            \midrule
                 0.15 & 1.0 & 0.25 & 1.0 & 2.0 & 0.5\\
            \bottomrule
            \end{tabular}
            \label{tab: coeffs 2}
        \end{table}

        \vspace{0.2cm}
        
        \begin{minipage}{\textwidth}
              \centering
              \begin{subfigure}{.5\textwidth}
              \centering
                  \includegraphics[width=.99\linewidth]{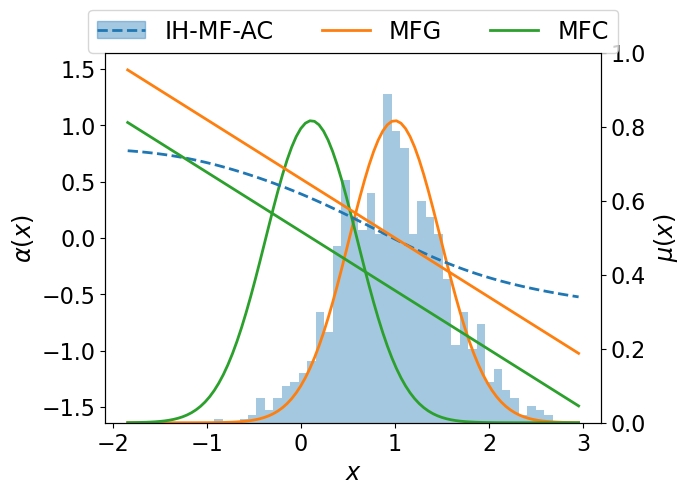}
                  \caption{MFG results for $\rho_\Sigma = 10^{-6}$}
              \label{fig: mfg results}
              \end{subfigure}%
            \begin{subfigure}{.5\textwidth}
              \centering
              \includegraphics[width=.99\linewidth]{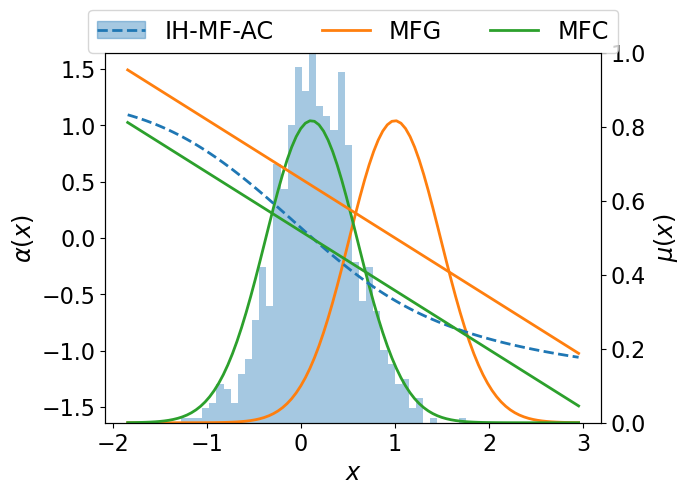}
              \caption{MFC results for $\rho_\Sigma = 5 \times 10^{-4}$} 
              \label{fig: mfc results}
              \end{subfigure}
              \caption{The histogram (blue) is the learned asymptotic distribution using samples generated from $\Sigma_{\varphi_n}$ and the dashed line (blue) is the learned feedback control after $N=10^6$ iterations averaged over five runs with different initial samples. The green curves correspond to the optimal control and mean field distribution for MFC, while the orange curves are the equivalent for MFG. The bottom axis shows the state variable $x$, the left axis refers to the value of the control $\alpha(x)$, and the right axis represents the probability density of $\mu$.}
              \label{fig: results 2}
        \end{minipage}

        \vspace{0.5cm}

        \begin{minipage}{\textwidth}
              \centering
              \begin{subfigure}{.5\textwidth}
              \centering
                  \includegraphics[width=.98\linewidth]{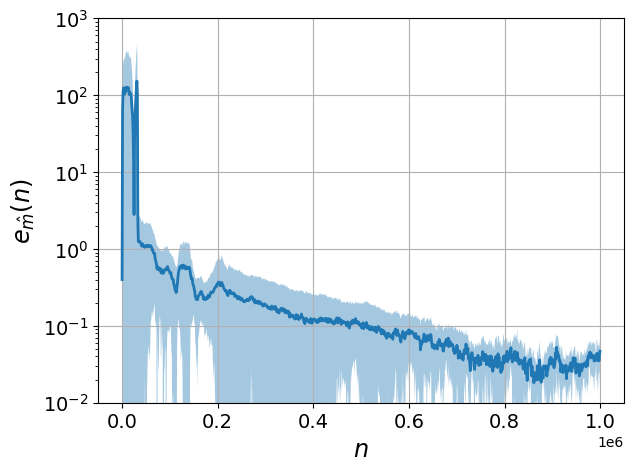}
                  \caption{Results for $\rho_\Sigma = 10^{-6}$}
              \label{fig: mfg mean error 2}
              \end{subfigure}%
            \begin{subfigure}{.5\textwidth}
              \centering
              \includegraphics[width=.98\linewidth]{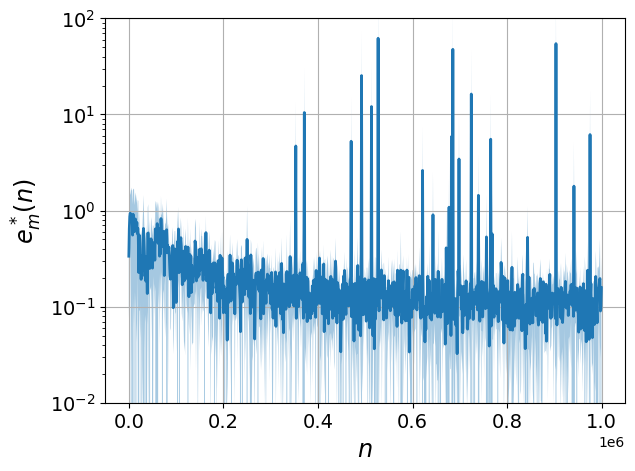}
              \caption{Results for $\rho_\Sigma = 5 \times 10^{-4}$}
              \label{fig: mfc mean error 2}
              \end{subfigure}
              \caption{We plot the absolute error between the mean of samples produced from the parameterized score function $\Sigma_{\varphi_n}$ and the optimal mean $\hat{m}$ in the case of MFG (left) and $m^*$ in the case of MFC (right). These values were averaged over five runs each with different random initial samples with the standard deviation given by the light blue shaded region. Large jumps are due to random outliers which result from the stochasticity of our algorithm.}
              \label{fig: mean errors 2}
        \end{minipage}
    \end{minipage}
\end{figure}

\begin{figure}[H]
    \centering
    \begin{minipage}{\textwidth}

        \begin{table}[H]
        \caption{Running cost coefficients and volatility for \cref{eq: multivar cost,eq: multivar dynamics}. The results for this parameter set are displayed in \cref{fig: 2d controls,fig: 2d dist,fig: 2d mean errors,fig: 2d control errors,fig: 2d value errors,fig: 2d values}. The coefficient values below were randomly generated with the exception of $\sigma$.}
        \centering
            \begin{tabular}{cccccc}
            \toprule
                 $C_1$ & $C_2$ & $C_3$ & $c_4$ & $C_5$ & $\sigma$ \\
            \midrule
                 $\begin{bmatrix}
                     0.964 & 0.236\\0.236 & 0.076
                 \end{bmatrix}$ &
                 $\begin{bmatrix}
                     0.677 & 0.937\\0.937 & 1.357
                 \end{bmatrix}$&
                 $\begin{bmatrix}
                     0.988 & 1.188\\1.188 & 1.483
                 \end{bmatrix}$ &
                 $\begin{bmatrix}
                     0.810\\ 0.872
                 \end{bmatrix}$ &
                 $\begin{bmatrix}
                     0.011 & 0.072\\0.072 & 0.520
                 \end{bmatrix}$&
                 $\begin{bmatrix}
                     0.3 & 0\\0 & 0.3
                 \end{bmatrix}$\\
            \bottomrule
            \end{tabular}
            \label{tab: coeffs 3}
        \end{table}
        
        \vspace{0.1cm}

        \begin{minipage}{\textwidth}
        \centering
            \begin{subfigure}{.5\textwidth}
            \centering
                \includegraphics[width=\linewidth]{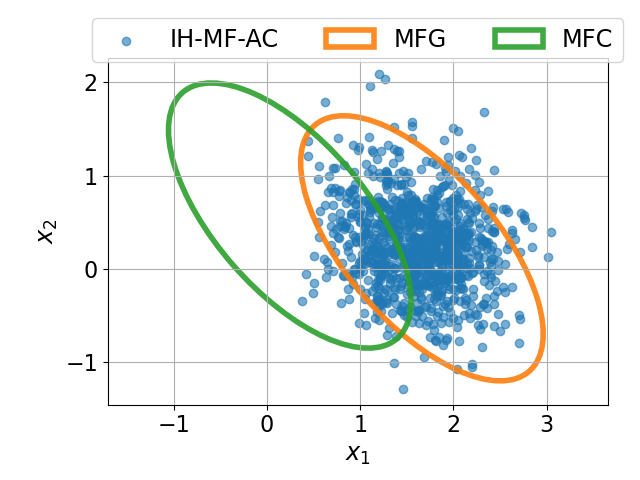}
                \vspace*{-0.8cm}
                \caption{Results for $\rho_\Sigma = 10^{-6}$}
                \label{fig: 2d mfg dist}
            \end{subfigure}%
            \begin{subfigure}{.5\textwidth}
                \centering
                  \includegraphics[width=\linewidth]{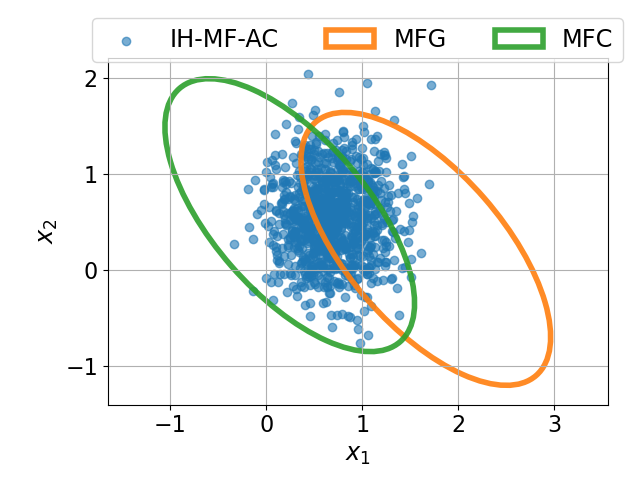}
                  \vspace*{-0.8cm}
                  \caption{Results for $\rho_\Sigma = 5 \times 10^{-4}$}
                  \label{fig: 2d mfc dist}
            \end{subfigure}
            \caption{The scatter plot of points (blue) are the samples of the learned asymptotic distribution generated from the parameterized score function $\Sigma_{\varphi_n}$ after $N=10^6$ iterations. The solid ellipses are the set of points with Mahalanobis distance three from the optimal mean field distributions in the case of MFG (orange) and MFC (green).}
            \label{fig: 2d dist}
        \end{minipage}

        \vspace{0.5cm}

        \begin{minipage}{\textwidth}
            \centering
            \begin{subfigure}{.5\textwidth}
                \centering
                \includegraphics[width=\linewidth]{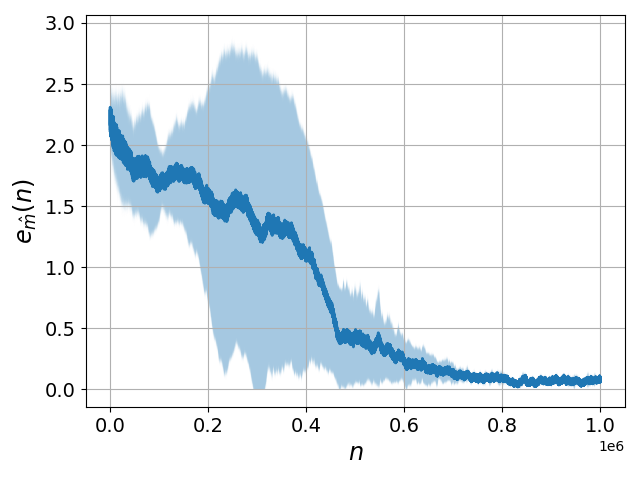}
                \vspace*{-0.8cm}
                \caption{Results for $\rho_\Sigma = 10^{-6}$}
                \label{fig: 2d mfg mean error}
            \end{subfigure}%
            \begin{subfigure}{.5\textwidth}
                \centering
                \includegraphics[width=\linewidth]{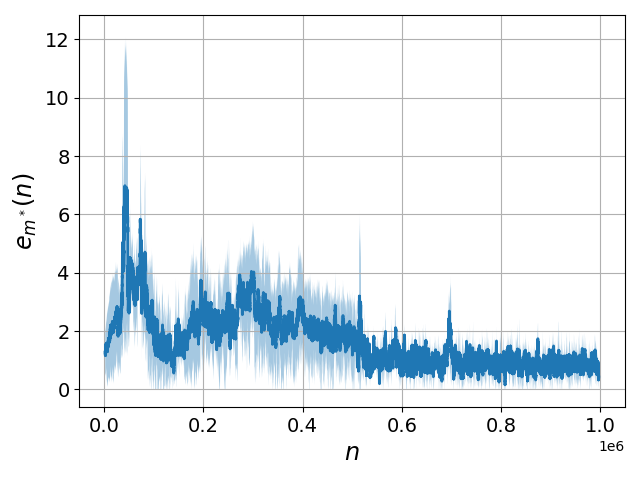}
                \vspace*{-0.8cm}
                \caption{Results for $\rho_\Sigma = 5 \times 10^{-4}$}
                \label{fig: 2d mfc mean error}
            \end{subfigure}
            \caption{We plot the absolute error between the mean of samples produced from the parameterized score function $\Sigma_{\varphi_n}$ and the optimal mean $\hat{m}$ in the case of MFG (left) and $m^*$ in the case of MFC (right). These values were averaged over five runs each with different random initial samples with the standard deviation given by the light blue shaded region. Large jumps are due to random outliers which result from the stochasticity of our algorithm.}
            \label{fig: 2d mean errors}
        \end{minipage}
        
    \end{minipage}
\end{figure}

For our numerical experiment, we test our algorithm on two different sets of values for the running cost coefficients and volatility $\sigma$ as listed in \Cref{tab: coeffs 1,tab: coeffs 2} for the 1-dimensional case and \Cref{tab: coeffs 3} for the 2-dimensional case. The discount factor is fixed in all cases to $\beta = 1$, and the continuous time is discretized using step size $\Delta t = 0.01$. The critic and score functions are both feedforward neural networks with one hidden layer of 128 neurons and a \texttt{tanh} activation function. For the 1-dimensional examples, the actor is also a feedforward neural network that outputs the mean and standard deviation of a normal distribution from which an action is sampled. Its architecture consists of a shared hidden layer of size 64 neurons and a \texttt{tanh} activation followed by two separate layers of size 64 neurons for the mean and standard deviation. The standard deviation layer is bookended by a \texttt{softmax} activation function to ensure its output is positive. The actor is meant to converge to a deterministic policy---also known as a pure control---over time, so in order to ensure a minimal level of exploration, we add a baseline value of $10^{-5}$ to the output layer. This straightforwardly mimics the notion of entropy regularization detailed in \cite{JMLR:v21:19-144}. The actor for the 2-dimensional example uses the same number of hidden layers and hidden neurons, and, correspondingly, outputs a mean vector and a diagonal covariance matrix to induce exploration. Refer to \Cref{tab: learning rates} for the learning rates used by the actor, critic, and score networks which were selected via grid search to give the best results. \Cref{tab: architectures} summarizes the total parameter count for each neural network.

For the Langevin Monte Carlo iterations, we pick a step size $\epsilon = 5 \times 10^{-2}$ as shown in \Cref{tab: learning rates}, and we run 200 iterations at each step using $k=1000$ samples.

\begin{figure}
    \begin{minipage}{\textwidth}
              \centering
              \begin{subfigure}{.5\textwidth}
              \centering
                  \includegraphics[width=.98\linewidth]{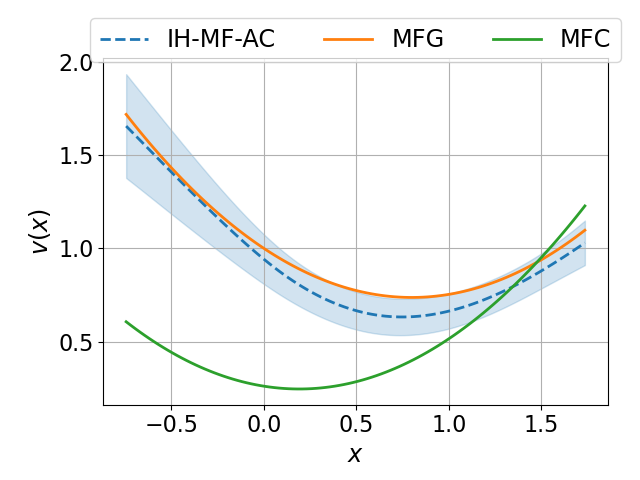}
                  \caption{Results for $\rho_\Sigma = 10^{-6}$}
              \label{fig: mfg value 1}
              \end{subfigure}%
            \begin{subfigure}{.5\textwidth}
              \centering
              \includegraphics[width=.98\linewidth]{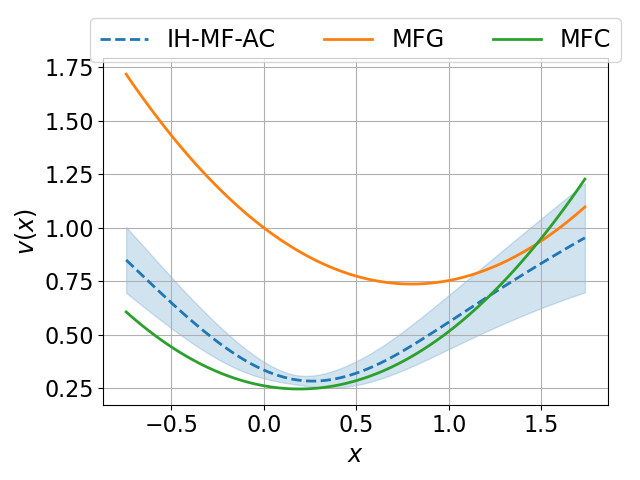}
              \caption{Results for $\rho_\Sigma = 5 \times 10^{-4}$}
              \label{fig: mfc value 1}
              \end{subfigure}
              \caption{The orange and green curves are the optimal value functions for the MFG and MFC problem, respectively. The blue dashed line is the learned value function given by the negative of critic $V_{\theta_N}$ averaged over five runs with different initial samples after $N = 10^6$ iterations. Since the original optimization problem aims to minimize cost while our algorithm seeks to maximize reward, we take the negative of the critic function to make the problems equivalent. The light blue shaded region depicts one standard deviation from the learned value.}
              \label{fig: values 1}
    \end{minipage}
    
    \vspace{0.8cm}
    
    \begin{minipage}{\textwidth}
              \centering
              \begin{subfigure}{.5\textwidth}
              \centering
                  \includegraphics[width=.98\linewidth]{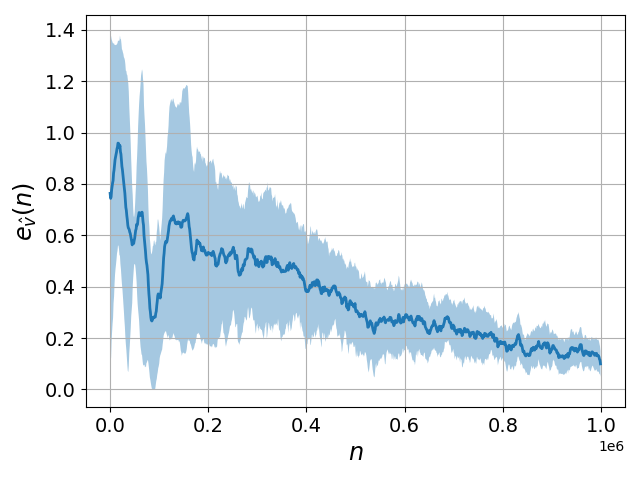}
                  \caption{Results for $\rho_\Sigma = 10^{-6}$}
              \label{fig: mfg value error 1}
              \end{subfigure}%
            \begin{subfigure}{.5\textwidth}
              \centering
              \includegraphics[width=.98\linewidth]{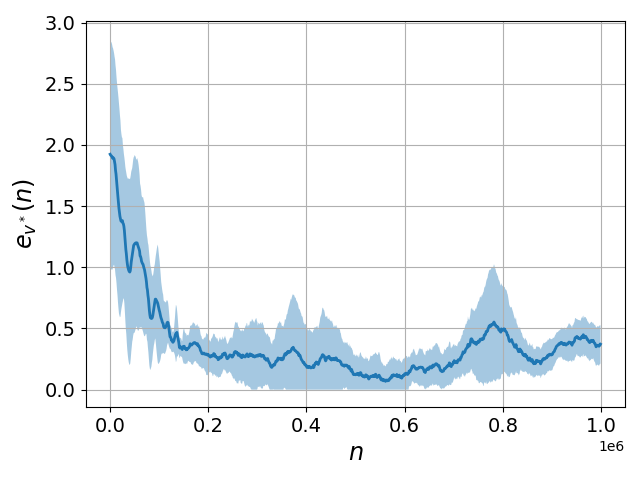}
              \caption{Results for $\rho_\Sigma = 5 \times 10^{-4}$}
              \label{fig: mfc value error 1}
              \end{subfigure}
              \caption{We plot the absolute error between the expected error between the learned value function $V_{\theta_n}$ and the optimal value function $\hat{v}$ in the case of MFG (left) and $v^*$ in the case of MFC (right). These plots were averaged over five runs each with different random initial samples and with the standard deviation given by the light blue shaded region. Large jumps are due to random outliers which result from the stochasticity of our algorithm.}
              \label{fig: value errors 1}
    \end{minipage}
\end{figure}

\begin{figure}[H]
    \begin{minipage}{\textwidth}
              \centering
              \begin{subfigure}{.5\textwidth}
              \centering
                  \includegraphics[width=.98\linewidth]{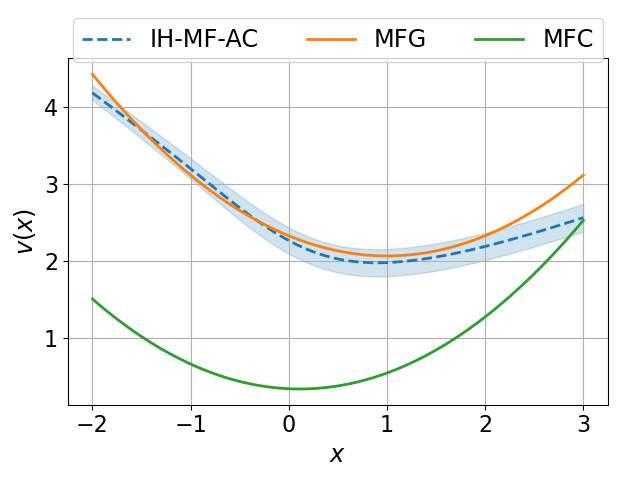}
                  \caption{Results for $\rho_\Sigma = 10^{-6}$}
              \label{fig: mfg value 2}
              \end{subfigure}%
            \begin{subfigure}{.5\textwidth}
              \centering
              \includegraphics[width=.98\linewidth]{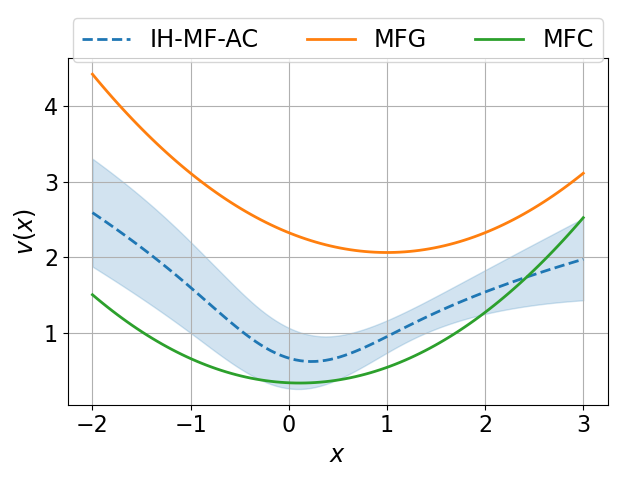}
              \caption{Results for $\rho_\Sigma = 5 \times 10^{-4}$}
              \label{fig: mfc value 2}
              \end{subfigure}
              \caption{The orange and green curves are the optimal value functions for the MFG and MFC problem, respectively. The blue dashed line is the learned value function given by the negative of critic $V_{\theta_N}$ averaged over five runs with different initial samples after $N = 10^6$ iterations. Since the original optimization problem aims to minimize cost while our algorithm seeks to maximize reward, we take the negative of the critic function to make the problems equivalent. The light blue shaded region depicts one standard deviation from the learned value.}
              \label{fig: values 2}
    \end{minipage}

    \vspace{1.2cm}
    
    \begin{minipage}{\textwidth}
              \centering
              \begin{subfigure}{.5\textwidth}
              \centering
                  \includegraphics[width=.98\linewidth]{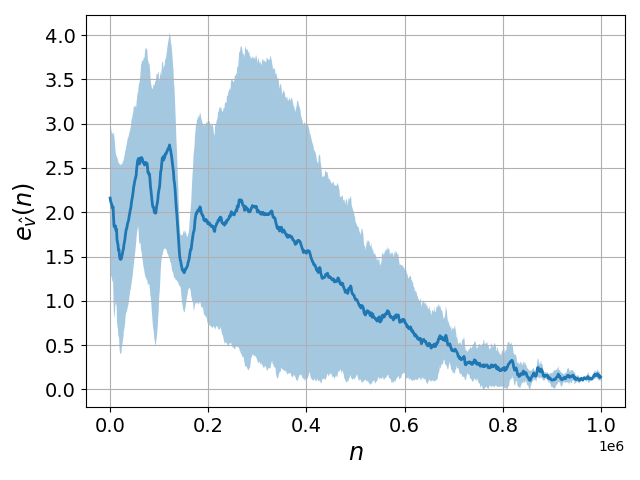}
                  \caption{Results for $\rho_\Sigma = 10^{-6}$}
              \label{fig: mfg value error 2}
              \end{subfigure}%
            \begin{subfigure}{.5\textwidth}
              \centering
              \includegraphics[width=.98\linewidth]{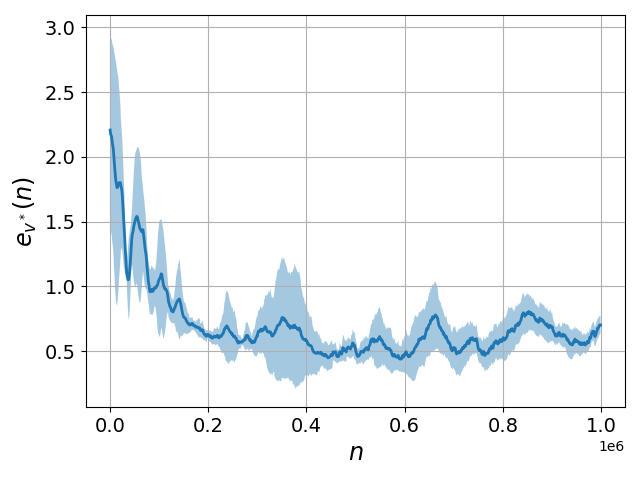}
              \caption{Results for $\rho_\Sigma = 5 \times 10^{-4}$}
              \label{fig: mfc value error 2}
              \end{subfigure}
              \caption{We plot the absolute error between the expected error between the learned value function $V_{\theta_n}$ and the optimal value function $\hat{v}$ in the case of MFG (left) and $v^*$ in the case of MFC (right). These plots were averaged over five runs each with different random initial samples and with the standard deviation given by the light blue shaded region. Large jumps are due to random outliers which result from the stochasticity of our algorithm.}
              \label{fig: value errors 2}
    \end{minipage}
\end{figure}

\begin{figure}[H]
    \centering
    \begin{minipage}{\textwidth}
        
        \begin{minipage}{\textwidth}
            \centering
            \begin{subfigure}{\textwidth}
                \includegraphics[width=\linewidth]{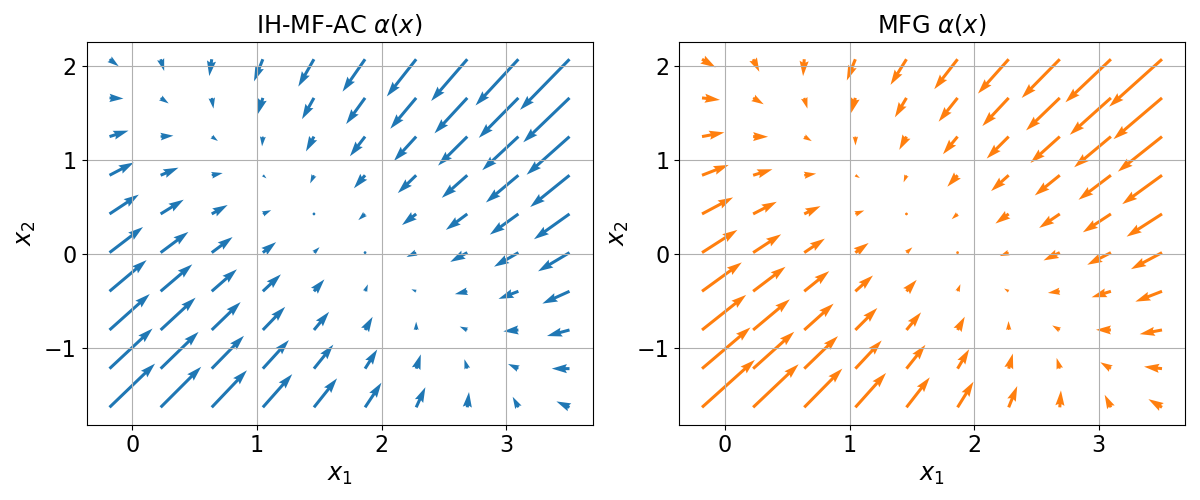}
                \caption{Results for $\rho_\Sigma = 10^{-6}$}
                \label{fig: 2d mfg control}
            \end{subfigure}\vspace{1cm}
            \begin{subfigure}{\textwidth}
                \includegraphics[width=\linewidth]{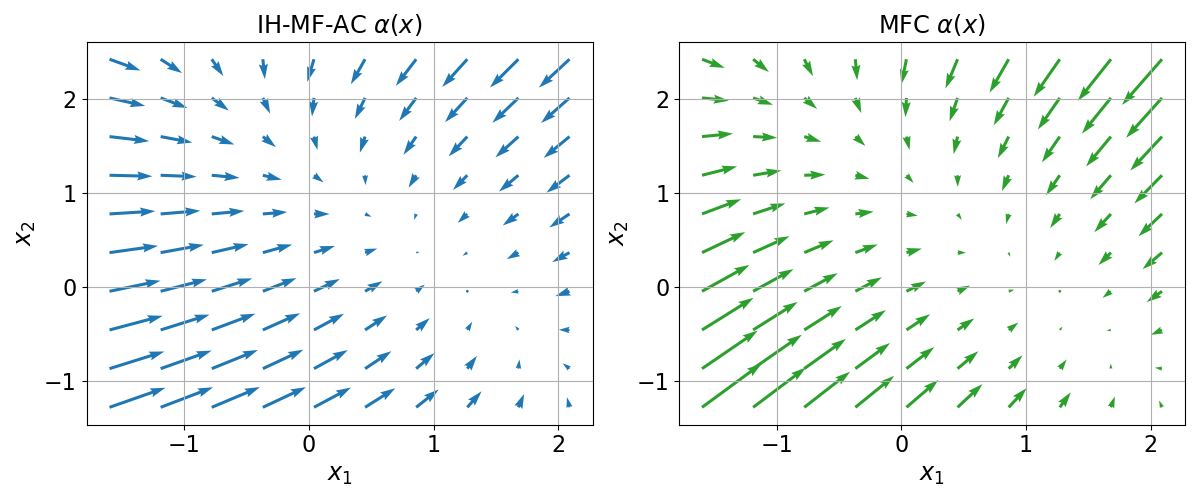}
                \caption{Results for $\rho_\Sigma = 5 \times 10^{-4}$}
                \label{fig: 2d mfc control}
            \end{subfigure}
            \caption{The orange and green vector fields (right) are the optimal controls for the MFG (top) and MFC (bottom) problems, respectively. The blue vector fields (left) show the learned feedback controls after $N = 10^6$ iterations averaged over five runs with different initial samples.}
            \label{fig: 2d controls}
        \end{minipage}
    \end{minipage}
\end{figure}

\begin{figure}[H]
    \centering
    \begin{minipage}{\textwidth}
        
        \begin{minipage}{\textwidth}
            \centering
            \begin{subfigure}{\textwidth}
                \includegraphics[width=\linewidth]{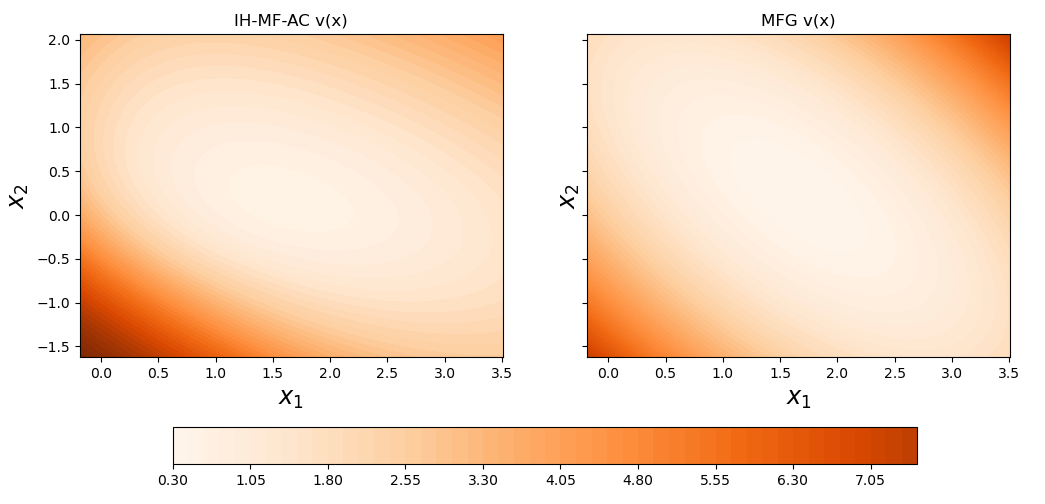}
                \caption{Results for $\rho_\Sigma = 10^{-6}$}
                \label{fig: 2d mfg value}
            \end{subfigure}\vspace{1cm}
            \begin{subfigure}{\textwidth}
                \includegraphics[width=\linewidth]{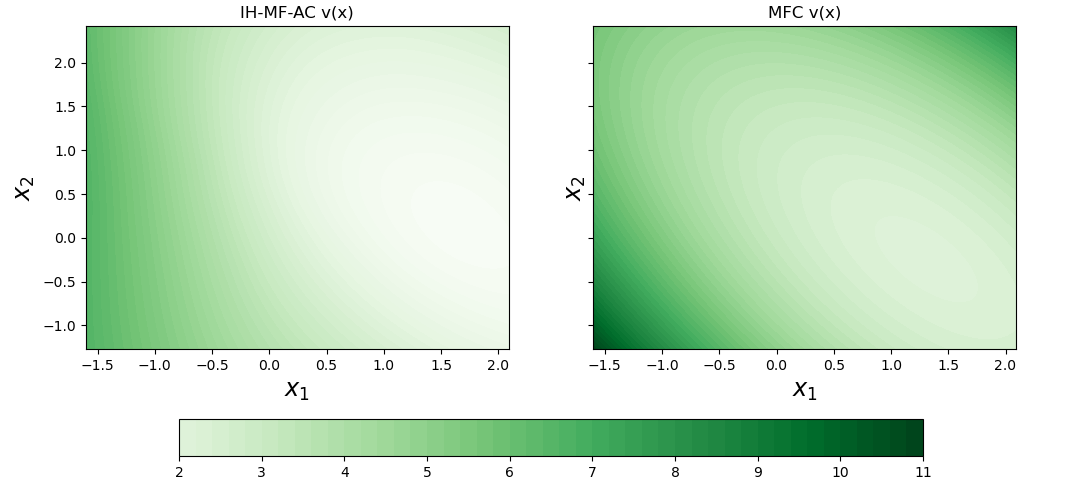}
                \caption{Results for $\rho_\Sigma = 5 \times 10^{-4}$}
                \label{fig: 2d mfc value}
            \end{subfigure}
            \caption{The right-hand side plots are the optimal value functions for the MFG (top) and MFC (bottom) problems, respectively. The left-hand side plots show the learned functions after $N = 10^6$ iterations averaged over five runs with different initial samples.}
            \label{fig: 2d values}
        \end{minipage}
    \end{minipage}
\end{figure}

The results of the algorithm applied to the LQ benchmark problem after $N=10^6$ iterations are displayed in \cref{fig: results 1,fig: values 1,fig: results 2,fig: values 2,fig: 2d dist,fig: 2d controls,fig: 2d values} with different sets of parameters along with the corresponding analytic solutions.
In addition to the plots of the learned solutions, we also provide plots for various error values to observe convergence with respect to the number of training steps. We consider the following error metrics for the MFG problem and their corresponding counterparts for the MFC problem:

\begin{itemize}
    \item The absolute error in the learned population mean:
    \[
        e_{\hat{m}}(n) \coloneqq ||m_{t_n} - \hat{m}||_2,
    \]

    \item The expected absolute error in the learned control:
    \[
        e_{\hat{\alpha}}(n) \coloneqq \E_{x \sim \hat{\mu}}\left[ ||\alpha_{\psi_n}(x) - \hat{\alpha}(x)||_2 \right] \approx \frac{1}{\ell} \sum_{j = 1}^\ell ||\alpha_{\psi_n}(x_j) - \hat{\alpha}(x_j)||_2,
    \]

    \item The expected absolute error in the learned value function:
    \[
        e_{\hat{v}}(n) \coloneqq \E_{x \sim \hat{\mu}}\left[ ||V_{\theta_n}(x) - \hat{v}(x)||_2 \right] \approx \frac{1}{\ell} \sum_{j = 1}^\ell ||V_{\theta_n}(x_j) - \hat{v}(x_j)||_2,
    \]
\end{itemize}
where $x_j$ are i.i.d. with distribution $\hat{\mu}$ for $j=1, 2, \dots, \ell$.

\begin{table}[tb]
\caption{Choice of learning rates used to obtain results in all figures seen in this work---$\rho_\Pi$ for the actor,  $\rho_V$ for the critic, $\rho_\Sigma$ for the score function, and $\epsilon$ for the Langevin dynamics. The two left columns are the learning rates used for the case $d=1$ and the right two for $d=2$. Boldface values indicate a difference in the learning rate between the MFG and MFC regimes.}
    \centering
    \begin{tabular}{lcc|cc}
        \toprule
        &$d=1$ & & $d=2$ &\\
        \midrule
         & MFG & MFC & MFG & MFC \\
        \midrule
        $\rho_\Pi$ & $5 \times 10^{-6}$ & $5 \times 10^{-6}$ & $5 \times 10^{-6}$ & $5 \times 10^{-6}$\\
        $\rho_V$ & $10^{-5}$ & $10^{-5}$ & $10^{-5}$ & $10^{-5}$ \\
        $\rho_{\Sigma}$ & $\mathbf{10^{-6}}$ & $\mathbf{5 \times 10^{-4}}$ & $\mathbf{10^{-6}}$ & $\mathbf{10^{-5}}$\\
        $\epsilon$ & $5\times 10^{-2}$ & $5\times 10^{-2}$ & $5\times 10^{-2}$ & $5\times 10^{-2}$\\
        \bottomrule
    \end{tabular}
    \label{tab: learning rates}
\end{table}

\begin{table}[tb]
\caption{Parameter counts and activation functions for the actor $\Pi_{\psi}$, critic $V_{\theta}$, and score $\Sigma_{\varphi}$ neural networks used to obtain results in all of the figures seen in this work.}
    \centering
    \begin{tabular}{lccc}
         \toprule
         & Actor & Critic & Score\\
         \midrule
         \# parameters & 258 & 385 & 385\\
         activation & \texttt{tanh} & \texttt{ELU} & \texttt{tanh}\\
         \bottomrule
    \end{tabular}
    \label{tab: architectures}
\end{table}\textbf{}

It should be noted that taking larger values of the volatility $\sigma$ yields degrading results in the learning of the solution as a result of the lower signal-to-noise ratio in the state samples. \Cref{tab: sigma mfg,tab: sigma mfc} illustrate the increase in several error metrics at the final time step $N$ for the 1-dimensional experiment with all other values fixed.

\begin{figure}
    \begin{minipage}{\textwidth}
              \centering
              \begin{subfigure}{.5\textwidth}
              \centering
                  \includegraphics[width=.98\linewidth]{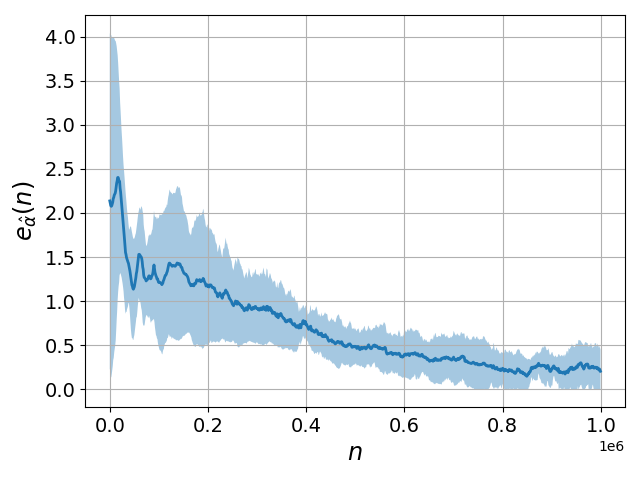}
                  \caption{Results for $\rho_\Sigma = 10^{-6}$}
              \label{fig: mfg control error 1}
              \end{subfigure}%
            \begin{subfigure}{.5\textwidth}
              \centering
              \includegraphics[width=.98\linewidth]{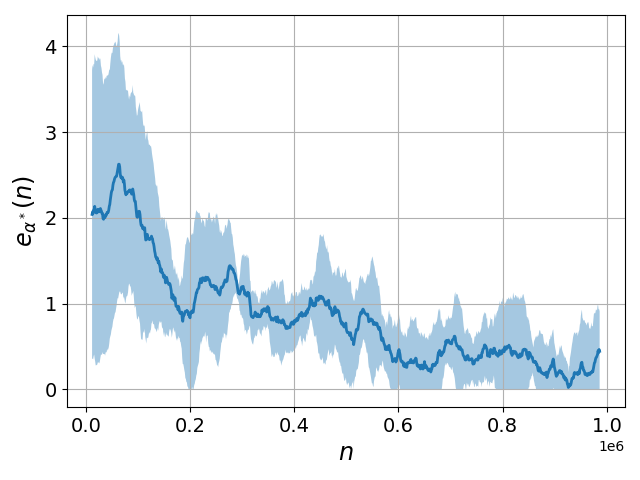}
              \caption{Results for $\rho_\Sigma = 5 \times 10^{-4}$}
              \label{fig: mfc control error 1}
              \end{subfigure}
              \caption{We plot the expected error between the learned control function $\alpha_n = \mathbb{E}[\Pi_{\psi_n}]$ and the optimal control $\hat{\alpha}$ in the case of MFG (left) and $\alpha^*$ in the case of MFC (right). These plots were averaged over five runs each with different random initial samples and with the standard deviation given by the light blue shaded region. Large jumps are due to random outliers which result from the stochasticity of our algorithm.}
              \label{fig: control errors 1}
    \end{minipage}
\end{figure}

\begin{figure}
    \begin{minipage}{\textwidth}
              \centering
              \begin{subfigure}{.5\textwidth}
              \centering
                  \includegraphics[width=.98\linewidth]{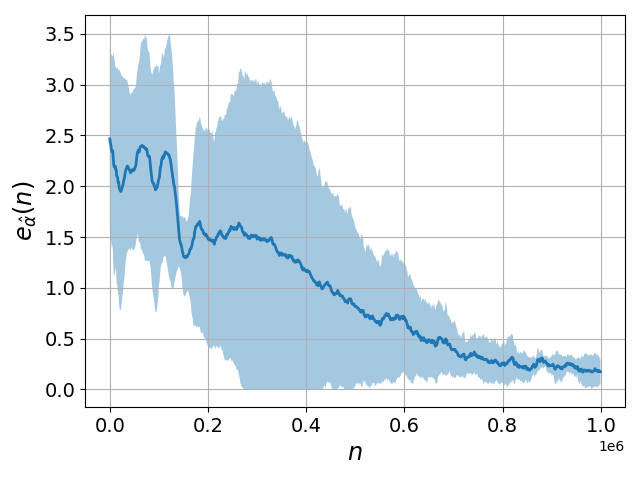}
                  \caption{Results for $\rho_\Sigma = 10^{-6}$}
              \label{fig: mfg control error 2}
              \end{subfigure}%
            \begin{subfigure}{.5\textwidth}
              \centering
              \includegraphics[width=.98\linewidth]{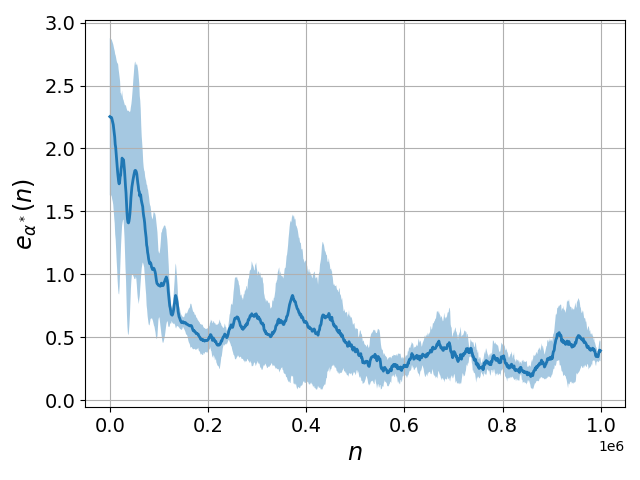}
              \caption{Results for $\rho_\Sigma = 5 \times 10^{-4}$}
              \label{fig: mfc control error 2}
              \end{subfigure}
              \caption{We plot the expected error between the learned control function $\alpha_n = \mathbb{E}[\Pi_{\psi_n}]$ and the optimal control $\hat{\alpha}$ in the case of MFG (left) and $\alpha^*$ in the case of MFC (right). These plots were averaged over five runs each with different random initial samples and with the standard deviation given by the light blue shaded region. Large jumps are due to random outliers which result from the stochasticity of our algorithm.}
              \label{fig: control errors 2}
    \end{minipage}
\end{figure}

\begin{figure}[H]
    \centering
    \begin{minipage}{\textwidth}
        \begin{minipage}{\textwidth}
            \centering
            \begin{subfigure}{.5\textwidth}
                \centering
                \includegraphics[width=\linewidth]{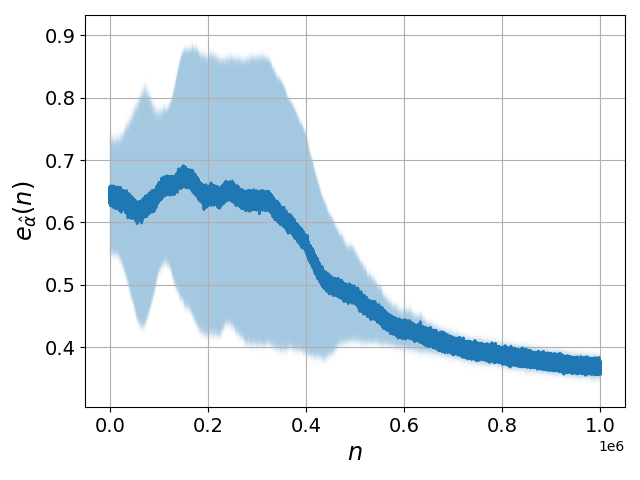}
                \caption{Results for $\rho_\Sigma = 10^{-6}$}
                \label{fig: 2d mfg control error}
            \end{subfigure}%
            \begin{subfigure}{.5\textwidth}
                \centering
                \includegraphics[width=\linewidth]{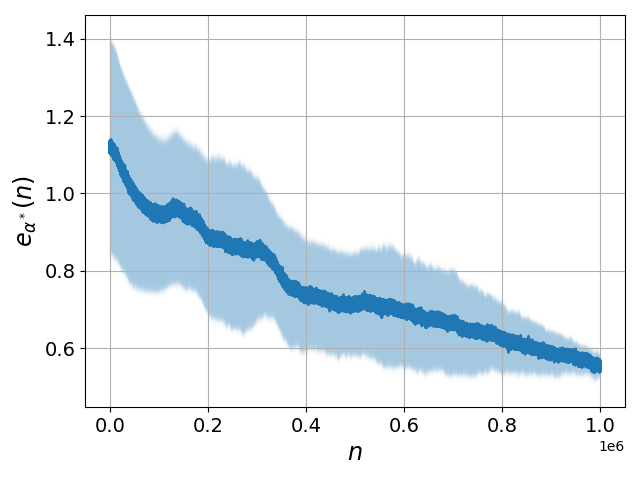}
                \caption{Results for $\rho_\Sigma = 10^{-5}$}
                \label{fig: 2d mfc control error}
            \end{subfigure}
            \caption{We plot the expected error between the learned control function $\alpha_n = \mathbb{E}[\Pi_{\psi_n}]$ and the optimal control $\hat{\alpha}$ in the case of MFG (left) and $\alpha^*$ in the case of MFC (right). These plots were averaged over five runs each with different random initial samples and with the standard deviation given by the light blue shaded region. Large jumps are due to random outliers which result from the stochasticity of our algorithm.}
            \label{fig: 2d control errors}
        \end{minipage}\vspace{1cm}
        \begin{minipage}{\textwidth}
            \centering
            \begin{subfigure}{.5\textwidth}
                \centering
                \includegraphics[width=\linewidth]{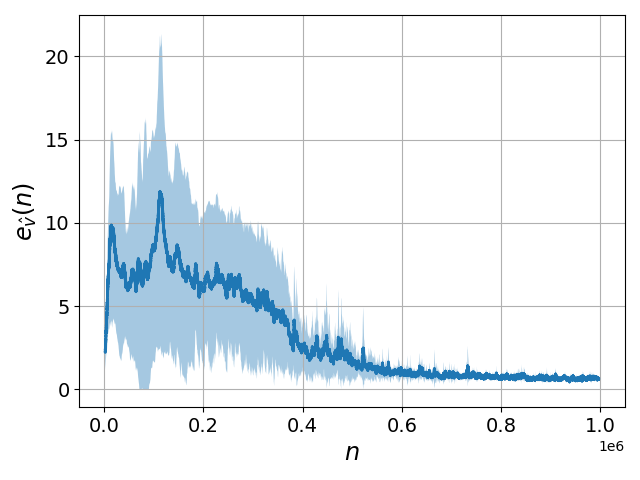}
                \caption{Results for $\rho_\Sigma = 10^{-6}$}
                \label{fig: 2d mfg value error}
            \end{subfigure}%
            \begin{subfigure}{.5\textwidth}
                \centering
                \includegraphics[width=\linewidth]{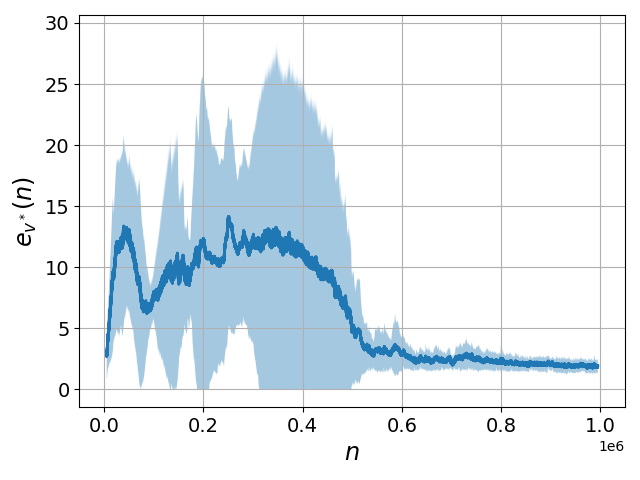}
                \caption{Results for $\rho_\Sigma = 10^{-5}$}
                \label{fig: 2d mfc value error}
            \end{subfigure}
            \caption{We plot the absolute error between the expected error between the learned value function $V_{\theta_n}$ and the optimal value function $\hat{v}$ in the case of MFG (left) and $v^*$ in the case of MFC (right). These plots were averaged over five runs each with different random initial samples and with the standard deviation given by the light blue shaded region. Large jumps are due to random outliers which result from the stochasticity of our algorithm.}
            \label{fig: 2d value errors}
        \end{minipage}
        
    \end{minipage}
\end{figure}

We observe many of the same insights alluded to by \cite{capponi_lehalle_2023,unified_q_learning} regarding the differences in recovering the MFG versus the MFC solution. Specifically, convergence to the MFG solution is more stable and faster than convergence to the MFC solution, as evidenced by the convergence plots in \cref{fig: mean errors 1,fig: mean errors 2,fig: 2d mean errors,fig: control errors 1,fig: control errors 2,fig: 2d control errors,fig: value errors 1,fig: value errors 2,fig: 2d value errors}. Further, in all cases, there were certain runs in which instability was amplified by the AC algorithm, in which case we saw the weights of the neural networks diverge to numerical overflow. In order to combat this, we imposed a bound on the state space during the first 200,000 iterations, truncating all states to the interval $[-5, 5]$ for 1-dimension case and the box $[-5,5] \times [-5, 5]$ in the 2-dimensional case. We removed the artificial truncation following the initial iterations and were able to mitigate the instability issues leading to overflow.

Observe that the optimal control is particularly well-learned within the support of the learned distribution. We postulate that a more intricate exploration scheme, perhaps along the lines of entropy regularization \cite{JMLR:v21:19-144}, may aid in learning the control in a larger domain. We conclude by noting that for all the numerical results in this paper, the gradient descent updates of \Cref{algo: ihmfac} (steps 5, 13, and 15) were computed using the Adam optimization update \cite{kingma2015} rather than the stochastic gradient descent update suggested in the pseudocode.
\newpage

 \clearpage

 \clearpage

\clearpage

\begin{table}[!tb]
    \parbox{0.47\linewidth}{
        \centering
        \begin{tabular}{c|c c c c c}
        \toprule
        $\sigma$ & 0.3  & 0.5 & 0.7 & 0.9 & 1.1\\
        \midrule
        $e_{\hat{m}}$ & 0.087 & 0.073 & 0.284 & 0.591 & 0.981 \\
        $e_{\hat{\alpha}}$ & 0.267 & 0.242 & 0.693 & 0.802 & 1.819 \\
        $e_{\hat{v}}$ & 0.178 & 0.394 & 0.522 & 0.673 & 2.005 \\
        \bottomrule
        \end{tabular}
        \caption{\revision{Error metrics for $\rho_{\Sigma} = 10^{-6}$. For increasing values of $\sigma$, we list the $L^2$ error in the learned mean $e_{\hat{m}}$, the expected $L^2$ error in the learned control $e_{\hat{\alpha}}$, and the expected $L^2$ error in the learned value function $e_{\hat{v}}$ after $N=10^6$ times steps. \label{tab: sigma mfg}}}
    }
    \hfill
    \parbox{0.48\linewidth}{
        \centering
        \begin{tabular}{c|c c c c c}
        \toprule
        $\sigma$ & 0.3  & 0.5 & 0.7 & 0.9 & 1.1\\
        \midrule
        $e_{m^*}$ & 0.115 & 0.299 & 1.876 & 2.701 & 3.415 \\
        $e_{\alpha^*}$ & 0.583 & 0.441 & 1.106 & 4.556 & 6.946 \\
        $e_{v^*}$ & 0.455 & 0.742 & 2.342 & 6.035 & 5.569 \\
        \bottomrule
        \end{tabular}
        \caption{\revision{Error metrics for $\rho_{\Sigma} = 5 \times 10^{-4}$. For increasing values of $\sigma$, we list the $L^2$ error in the learned mean $e_{m^*}$, the expected $L^2$ error in the learned control $e_{\alpha^*}$, and the expected $L^2$ error in the learned value function $e_{v^*}$ after $N=10^6$ time steps. \label{tab: sigma mfc}}}
    }
\end{table}

We have shown that the same algorithm tested on a continuous-state space infinite horizon LQ problem, with hyperparameters shown in \Cref{tab: learning rates}, recovers the MFG or MFC solution. However, our numerical experiments show that there is more stability for the MFG problem, as also observed in the case of the Q-learning algorithm in the context of infinite horizon and finite-space \cite{unified_q_learning}, and in \cite{capponi_lehalle_2023} in the context of discrete-space finite horizon problems. In the following section, we see that we observe a similar phenomenon in the case of the mixed mean field control game. Note, however, that in our general algorithm we are not taking advantage of the fact that the numerical example is linear-quadratic. Doing so, we would know a priori that the value function is quadratic and that the control is linear, as is done, for instance, in \cite{chen2023global}.


\section{Actor-Critic Algorithm for Mean Field Control Games (MFCG)} \label{sec: mfcg}
As observed in \cite{angiuli2023} in the case of tabular Q-learning, our IH-MF-AC algorithm (\Cref{algo: ihmfac}) can easily be extended to the case of mixed mean field control game problems that involve two population distributions, a local one and a global one. This type of game corresponds to competitive games among a large number of large groups, where agents within each group collaborate. The local distribution represents the distribution within each "representative" agent's group, while the global distribution represents the distribution of the entire population.

Such games are motivated as follows: Consider a finite population of agents consisting of $M$ groups, each containing $N$ agents. An agent indexed by $(m,n)$ indicates that she is the $n^{th}$ member of the $m^{th}$ group. Agents collaborate within their respective groups (sharing the same first index) and compete with all agents from other groups. In other words, all $N$ agents in group $m$ collectively work to minimize the total cost of group $m$. The MFCG corresponds to the mean-field limit as $N$ and $M$ tend to infinity, and serves as a proxy for the solution in the finite player case. We refer to  \cite{angiuli2023,angiuli2022c} for further details on MFCG, including the limit from finite player games to infinite player games. Note that the solution gives an approximation of the Nash equilibrium between the competitive groups.

The solution of an infinite horizon mean field control game is a control-mean field pair $(\hat{\alpha}, \hat{\mu}) \in \mathbb{A} \times  \mathcal{P}(\mathbb{R}^{d})$
satisfying the following:
\begin{enumerate}
    \item $\hat{\alpha}$ solves the McKean-Vlasov stochastic optimal control problem
    \begin{equation} \label{eq: mfcg cost}
        \inf _{\alpha \in \mathbb{A}} J_{\hat{\mu}}(\alpha)=\inf _{\alpha \in \mathbb{A}} \mathbb{E}\left[ \int_0^\infty e^{-\beta t} f\left(X_t^{\alpha, \hat{\mu}}, \hat{\mu},\mu^{\alpha, \hat{\mu}}, \alpha(X_t^{\alpha, \hat{\mu}})\right)\, \dd t \right],\quad \beta>0,
    \end{equation}
    subject to
    \begin{equation} \label{eq: mfcg dynamics}
        \dd X_t^{\alpha, \hat{\mu}}=b\left(X_t^{\alpha, \hat{\mu}}, \hat{\mu}, \mu^{\alpha, \hat{\mu}}, \alpha(X_t^{\alpha, \hat{\mu}})\right) \, \dd t+\sigma\left(X_t^{\alpha, \hat{\mu}}, \hat{\mu}, \mu^{\alpha, \hat{\mu}}, \alpha(X_t^{\alpha, \hat{\mu}})\right)\, \dd W_t, \quad X_0^{\alpha, \hat{\mu}}=\xi,
    \end{equation}
    where $\mu^{\alpha, \hat{\mu}}=\lim _{t \rightarrow \infty} \mathcal{L}(\mathrm{X}_t^{\alpha, \hat{\mu}})$;

    \item fixed point condition: $\hat{\mu} = \lim_{t \to \infty} \mathcal{L}(X_t^{\hat{\alpha}, \hat{\mu}})$.
\end{enumerate}
Note that conditions 1 and 2 above imply that $\hat{\mu} = \mu^{\hat{\alpha}, \hat{\mu}}$.

We modify \Cref{algo: ihmfac} into our \emph{infinite horizon mean field control game actor-critic} (IH-MFCG-AC) algorithm such that the global score function $\Sigma_\varphi$ represents the global distribution $\hat{\mu}$ and the local score function $\widetilde{\Sigma}_\xi$ represents the local distribution $\mu^{\alpha, \hat{\mu}}$. This is meant to mimic the parallel between the mean field game solution with the global distribution, and the mean field control solution with the local distribution. Following our intuition from \Cref{subsec: unifying}, our choice of the now four learning rates will be chosen according to
\begin{equation}\label{eq: mfcg lr relations}
    \rho_{\Sigma} < \min\{\rho_\Pi, \rho_V\} < \max\{\rho_\Pi, \rho_V\} < \rho_{\widetilde{\Sigma}}.
\end{equation} 
Refer to \Cref{algo: ihmfcgac} for the complete pseudocode.

\begin{algorithm}[!tb]
   \caption{\textbf{IH-MFCG-AC: Infinite Horizon Mean Field Control Game Actor-Critic}\\Differences from \Cref{algo: ihmfac} are highlighted in blue.}
   \label{algo: ihmfcgac}
\begin{algorithmic}[1] 
    \REQUIRE Number of time steps $N \gg 0$; discrete time step size $\Delta t$; neural network learning rates for actor $\rho_\Pi$, critic $\rho_V$, global score $\rho_\Sigma$, and local score $\rho_{\widetilde{\Sigma}}$; Langevin dynamics step size $\epsilon$.
    \STATE Initialize neural networks:\\
    \textbf{Actor} $\Pi_{\psi_0}: \R^d \to \mathcal{P}(\R^k)$\\
    \textbf{Critic} $V_{\theta_0}: \R^d \to \R$\\
    \textbf{Global Score} $\Sigma_{\varphi_0}: \R^d \to \R^d$\\
    {\color{blue}\textbf{Local Score} $\widetilde{\Sigma}_{\xi_0}: \R^d \to \R^d$}\vspace{0.2cm}
    \STATE Agent receives initial state $X_{t_0}$ from the Environment.\vspace{0.2cm}
   \FOR{$n=0,\dots,N-1$}\vspace{0.2cm}
      \STATE Environment computes score loss for $\Sigma$:
      $\quad L_\Sigma (\varphi_n) = \tr\left( \nabla_x \Sigma_{\varphi_n}(X_{t_n}) \right) + \frac{1}{2}\norm{\Sigma_{\varphi_n}(X_{t_n})}_2^2$\vspace{0.2cm}
      
      \STATE Environment updates $\Sigma$ with SGD:
      $\quad \varphi_{n+1} = \varphi_n -\rho_\Sigma \nabla_{\varphi} L_\Sigma (\varphi_n)$\vspace{0.2cm}

      {\color{blue}\STATE Environment computes score loss for $\widetilde{\Sigma}$:
      $\quad L_{\widetilde{\Sigma}} (\xi_n) = \tr\left( \nabla_x \widetilde{\Sigma}_{\xi_n}(X_{t_n}) \right) + \frac{1}{2}\norm{\widetilde{\Sigma}_{\xi_n}(X_{t_n})}_2^2$\vspace{0.2cm}

      \STATE Environment updates $\widetilde{\Sigma}$ with SGD:
      $\quad \xi_{n+1} = \xi_n -\rho_{\widetilde{\Sigma}} \nabla_{\xi} L_{\widetilde{\Sigma}} (\xi_n)$}\vspace{0.2cm}
      
      \STATE Environment generates mean field samples $S_{t_n} = \left(S_{t_n}^{(1)}, S_{t_n}^{(2)}, \dots, S_{t_n}^{(k)}\right)$ from $\Sigma_{\varphi_{n+1}}$ {\color{blue}and $\widetilde{S}_{t_n} = \left(\widetilde{S}_{t_n}^{(1)}, \widetilde{S}_{t_n}^{(2)}, \dots, \widetilde{S}_{t_n}^{(k)}\right)$ from $\widetilde{\Sigma}_{\xi_{n+1}}$} using Langevin dynamics (\cref{eq: langevin}) with step size $\epsilon$ and compute $\overline{\mu}_{S_{t_n}} \coloneqq \frac{1}{k} \sum_{i=1}^k \delta_{S_{t_n}^{(i)}}$ {\color{blue}and $\overline{\mu}_{\widetilde{S}_{t_n}} \coloneqq \frac{1}{k} \sum_{i=1}^k \delta_{\widetilde{S}_{t_n}^{(i)}}$}.\vspace{0.2cm}

      \STATE Agent samples action:
      $\quad A_{t_n} \sim \Pi_{\psi_n}(\cdot \mid X_{t_n})$\vspace{0.2cm}
      
      \STATE Agent observes reward $r_{n+1}$ and next state $X_{t_{n+1}}$from the environment.\vspace{0.2cm}

      \STATE Agent computes TD target:
      $\quad y_{n+1} = r_{n+1} + e^{-\beta \Delta t} V_{\theta_n}(X_{t_{n+1}})$\vspace{0.2cm}
      
      \STATE Agent computes TD error:
      $\quad \delta_{\theta_n} =y_{n+1} - V_{\theta_n}(X_{t_n})$\vspace{0.2cm}
      
      \STATE Agent computes critic loss:
      $\quad L_V(\theta_n) = \delta_{\theta_n}^2$\vspace{0.2cm}
      
      \STATE Agent updates critic with SGD:
      $\quad \theta_{n+1} = \theta_n - \rho_V \nabla_{\theta} L_V(\theta_n)$\vspace{0.2cm}
      
      \STATE Agent computes actor loss:
      $\quad L_{\Pi}(\psi_n) = -\delta_{\theta_n} \log \Pi_{\psi_n}(A_{t_n} \mid X_{t_n})$\vspace{0.2cm}
      
      \STATE Agent updates actor with SGD:
      $\quad \psi_{n+1} = \psi_n - \rho_{\Pi} \nabla_{\psi} L_\Pi(\psi_n)$\vspace{0.2cm}
   \ENDFOR \vspace{0.2cm}
   \RETURN $(\Pi_{\psi_N}, \Sigma_{\varphi_N}, \widetilde{\Sigma}_{\xi_N})$
\end{algorithmic}
\end{algorithm} \clearpage

\subsection{A Linear-Quadratic Benchmark}
We test \Cref{algo: ihmfcgac} on the following linear-quadratic MFCG. We wish to minimize
\begin{equation}\label{eq: mfcg lq cost}
    \begin{split}
        \mathbb{E}\Biggl[\int_0^{\infty} e^{-\beta t}\biggl(\frac{1}{2} \alpha_t^2+c_1\left(\mathrm{X}_t^{\alpha, \mu}-c_2 m\right)^2+c_3\left(\mathrm{X}_t^{\alpha, \mu}-c_4\right)^2\\
        \quad {}+\tilde{c}_1\left(\mathrm{X}_t^{\alpha, \mu}-\tilde{c}_2 m^{\alpha, \mu}\right)^2+\tilde{c}_5\left(m^{\alpha, \mu}\right)^2\biggr) \mathrm{d} t\Biggr]
    \end{split}
\end{equation}
subject to the dynamics
\begin{equation} \label{eq: mfcg lq dynamics}
    \mathrm{d}X_t^{\alpha, \mu}=\alpha_t \,\mathrm{d} t+\sigma \, \mathrm{d}W_t, \qquad t \in [0, \infty)
\end{equation}
where $m = \int x \, \dd \mu(x)$ and $m^{\alpha,\mu} = \int x \, \dd \mu^{\alpha, \mu}(x)$ and the fixed point condition $m=\lim _{t \rightarrow \infty} \mathbb{E}(X_t^{\hat{\alpha}, \mu})=m^{\hat{\alpha}, \mu}$ where $\hat{\alpha}$ is the optimal action.

We present the analytic solution to the MFCG problem using notation consistent with the derivation in \cite{angiuli2023}. The value function is defined as
\begin{equation}
    \begin{split}
        v(x) \coloneqq \inf _{\alpha \in \mathbb{A}} \mathbb{E}\Biggl[\int_0^{\infty} e^{-\beta t}\biggl(\frac{1}{2} \alpha_t^2+c_1\left(\mathrm{X}_t^{\alpha, \mu}-c_2 m\right)^2+c_3\left(\mathrm{X}_t^{\alpha, \mu}-c_4\right)^2\\
        \quad {}+\tilde{c}_1\left(\mathrm{X}_t^{\alpha, \mu}-\tilde{c}_2 m^{\alpha, \mu}\right)^2+\tilde{c}_5\left(m^{\alpha, \mu}\right)^2\biggr) \mathrm{d} t \mid X_0 = x\Biggr].
    \end{split}
\end{equation}
The explicit formula $v(x) = \Gamma_2 x^2 + \Gamma_1 x + \Gamma_0$ can be derived as the solution to the Hamilton-Jacobi-Bellman equation where
\begin{align*}
    \Gamma_2 &=\frac{-\beta+\sqrt{\beta^2+8\left(c_1+c_3+\tilde{c}_1\right)}}{4}\\[1em]
    \Gamma_1 &= -\frac{2 \Gamma_2 c_3 c_4}{c_1\left(1-c_2\right)+\tilde{c}_1\left(1-\tilde{c}_2\right)^2+c_3+\tilde{c}_5}\\[1em]
    \Gamma_0 &= \frac{c_1 c_2^2 m^2+\left(\tilde{c}_1 \tilde{c}_2^2+\tilde{c}_5\right)\left(m^{\alpha, \mu}\right)^2+\sigma^2 \Gamma_2-\frac{1}{2} \Gamma_1^2+c_3 c_4^2}{\beta} .
\end{align*}
Then the optimal control for the MFCG is
\begin{equation}\label{eq: mfcg optimal control}
    \hat{\alpha}(x) = -(2 \Gamma_2 x + \Gamma_1).
\end{equation}
Substituting \cref{eq: mfcg optimal control} into \cref{eq: mfcg lq dynamics} yields the Ornstein-Uhlenbeck process
\[
    \dd X_t = -\left(2 \Gamma_2 X_t + \Gamma_1\right)\, \dd t + \sigma \, \dd W_t
\]
whose limiting distribution is
\begin{equation}
    \hat{\mu} = \mu^{\hat{\alpha}, \hat{\mu}} = \mathcal{N}\left( -\frac{\Gamma_1}{2\Gamma_2}, \frac{\sigma^2}{4 \Gamma_2} \right).
\end{equation}
We note that an equation for $\hat{m}$ and $m^{\hat{\alpha}, \hat{\mu}}$ that only depends on the running cost coefficients is
\begin{equation} \label{eq: mfcg mean}
    m \coloneqq \hat{m} = m^{\hat{\alpha}, \hat{\mu}} = \frac{c_3 c_4}{c_1\left(1-c_2\right)+\tilde{c}_1\left(1-\tilde{c}_2\right)^2+c_3+\tilde{c}_5}.
\end{equation}

\subsection{Hyperparameters and Numerical Specifics}
For the LQ benchmark problem, we consider the following choice of parameters: $c_1 = 0.5$, $c_2 = 1.5$, $c_3 = 0.5$, $c_4 = 0.25$, $\tilde{c}_1 = 0.3$, $\tilde{c}_2 = 1.25$, $\tilde{c}_5 = 0.25$, discount factor $\beta = 1$, and volatility $\sigma = 0.5$. The time discretization is again $\Delta t = 0.01$. Our intention was to modify as few of the numerical hyperparameters from \Cref{sec: numerical} as possible, including the neural network architectures for the actor and critic. The global and local score networks both inherit the architecture from the score network described in \Cref{subsec: mf hyperperams} and \Cref{tab: architectures}. The learning rates for the networks are taken directly from \Cref{tab: learning rates} with the global and score network learning rates assuming the values used to obtain the MFG and MFC results, respectively, from \Cref{subsec: mf hyperperams}. This is to say, $(\rho_\Pi, \rho_V, \rho_\Sigma, \rho_{\widetilde{\Sigma}}) = (5\times 10^{-6}, 10^{-5}, 10^{-6}, 5 \times 10^{-4}$), which satisfy $\rho_\Sigma < \rho_\Pi < \rho_V < \rho_{\widetilde{\Sigma}}$, the learning rate inequality proposed in \cref{eq: mfcg lr relations}. The global and local distribution samples are computed at each time step using Langevin dynamics with $\epsilon = 5 \times 10^{-2}$ for 200 iterations using $k=1000$ samples.

\begin{figure}[!tb]
    \begin{minipage}{\textwidth}
        \centering
        \begin{subfigure}{\textwidth}
            \includegraphics[width=0.98\linewidth]{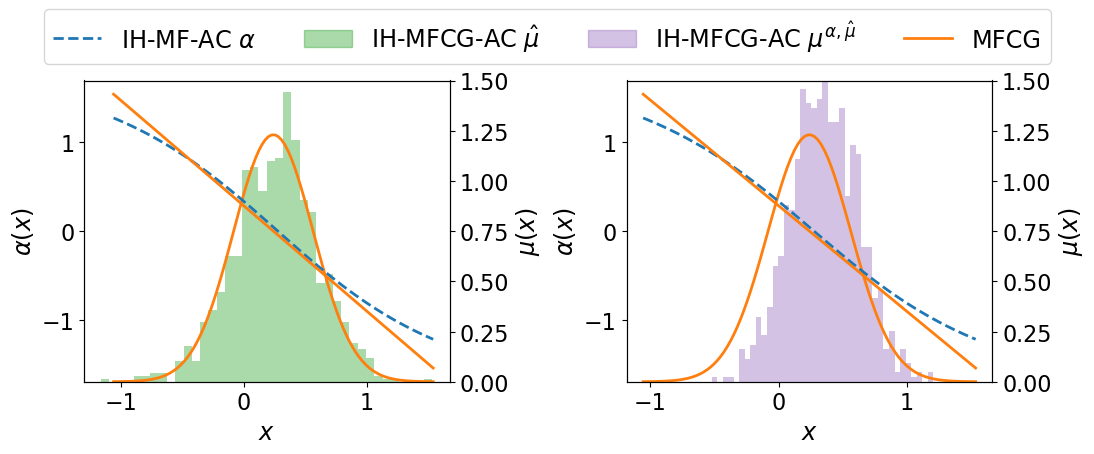}
        \end{subfigure}
        \caption{The histograms are the learned distributions generated using samples from the global score $\Sigma_{\varphi_n}$ (green) representing the global distribution $\hat{\mu}$ and the local score $\widetilde{\Sigma}_{\xi_n}$ (purple) representing the local score $\mu^{\alpha, \hat{\mu}}$ after $N = 2\times 10^6$ iterations. The dashed line (blue) is the learned feedback control averaged over five runs with different initial samples. The benchmark solution to the MFCG is provided in orange. The $x$-axis shows the state variable $x$, the left $y$-axis refers to the value of the control $\alpha(x)$, and the right axis represents the probability density of $\mu(x)$.}
        \label{fig: mfcg results}
        \end{minipage}
\end{figure}

\begin{figure}[!htb]
    \centering
    \begin{minipage}{\textwidth}

        \begin{minipage}{\textwidth}
        \centering
        \begin{subfigure}{\textwidth}
            \includegraphics[width=0.98\linewidth]{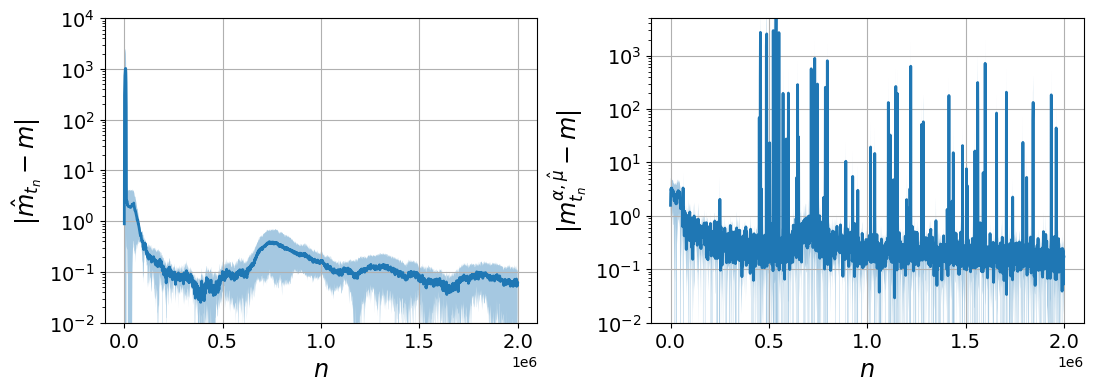}
        \end{subfigure}
        \caption{The blue curve is a rolling average of the absolute error of the mean of samples produced from the global score function $\Sigma_{\varphi_n}$ (left)---denoted $\hat{m}_{t_n}$---and the local score function $\widetilde{\Sigma}_{\xi_n}$---denoted $m^{\alpha, \hat{\mu}}_{t_n}$---compared to the optimal mean $m$ from \cref{eq: mfcg mean}. These values were averaged over five runs each with different random initial samples with the standard deviation given by the light blue shaded region. Large jumps are due to random outliers which result from the stochasticity of our algorithm.}
        \label{fig: mfcg means}
        \end{minipage}
        
        \vspace{2cm}

        \begin{minipage}{\textwidth}
              \centering
              \includegraphics[width=0.5\linewidth]{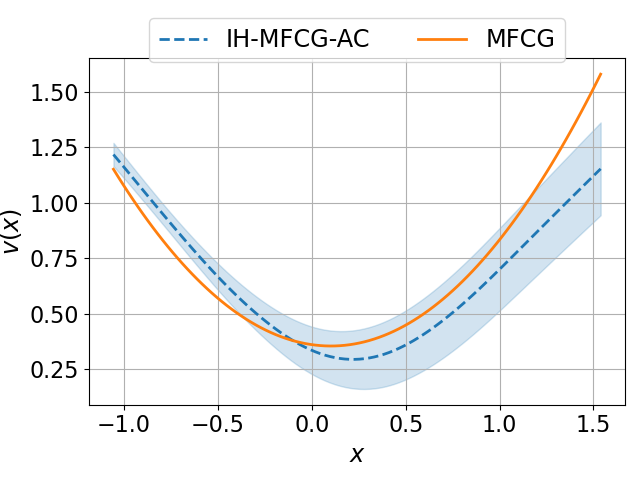}
              \caption{The orange curve is the optimal value function for the MFCG problem. The blue dashed line is the learned value function given by the negative of critic $V_{\theta_N}$ averaged over five runs with different initial samples after $N = 2 \times 10^6$ iterations. Since the original optimization problem aims to minimize cost while our algorithm seeks to maximize reward, we take the negative of the critic function to make the problems equivalent. The light blue shaded region depicts one standard deviation from the learned value.}
              \label{fig: mfcg values}
    \end{minipage}
    \end{minipage}
\end{figure}

The results of the IH-MFCG-AC algorithm (\Cref{algo: ihmfcgac}) are presented in \cref{fig: mfcg results,fig: mfcg means,fig: mfcg values}. As expected, the learning of the global and local distributions reflects that of the optimal MFG distribution and the optimal MFC distribution, respectively. We observe that the global score is learned faster and with more accuracy than the local score, which is prone to outliers and instability. The optimal control is learned well within the support of the optimal distribution, but could possibly be expanded with a more advanced exploration strategy.

 \clearpage

\section{Conclusion}
We have introduced a novel AC algorithm for solving infinite horizon mean field games and mean field control problems in continuous spaces. This algorithm, called IH-MF-AC, uses neural networks to parameterize a policy and value function, from which an optimal control is derived, as well as a score function, which represents the optimal mean field distribution on a continuous space. The MFG or MFC solution is arrived at depending on the choice of learning rates for the actor, critic, and score networks. We test our algorithm against a linear-quadratic benchmark problem and are able to recover the analytic solutions with a high degree of accuracy. Finally, we propose and test a modification of the algorithm, called IH-MFCG-AC, to solve the recently developed mixed mean field control game problems. For future work, several directions are worth investigating: enhancing the algorithm for scenarios where the signal-to-noise ratio is low (i.e., $\sigma$ is large) or where the horizon is finite, as well as conducting a rigorous numerical analysis of the proposed algorithm.

\section*{Acknowledgment}
J.F. was supported by NSF grant DMS-1953035. R.H. was partially supported by the NSF grant DMS-1953035, the Regents' Junior Faculty Fellowship at UCSB, and a grant from the Simons Foundation (MP-TSM-00002783). Use was made of computational facilities purchased with funds from the National Science Foundation (CNS-1725797) and administered by the Center for Scientific Computing (CSC). The CSC is supported by the California NanoSystems Institute and the Materials Research Science and Engineering Center (MRSEC; NSF DMR 1720256) at UC Santa Barbara. R.H. is grateful to Jingwei Hu for the useful discussions.

\bibliographystyle{apalike}

\bibliography{reference}
\newpage

\appendix
\section{Solution for Multivariate Linear-Quadratic Asymptotic MFG and MFC Problems} \label{ap: derivations}

For this problem, we assume the state process takes values in $\R^d$. Mimicking the notation established in \Cref{subsec: lq bench}, let $C_1 , C_3, C_5 \in \R^{d \times d}$ be symmetric positive-definite matrices, $c_4 \in \R^d$, and $C_2 \in \R^{d \times d}$. Further, let $W_t$ be an $m$-dimensional Brownian motion and $\sigma \in \R^{d \times m}$. For this class of mean field problems, the running cost is given by
\begin{equation} \label{eq: multivar mf cost}
    f(x, \mu, \alpha) = \frac{1}{2} \alpha^T \alpha + (x - C_2 m)^T C_1 (x - C_2 m) + (x - c_4)^T C_3 (x - c_4) + m^T C_5 m
\end{equation}
and the state dynamics by
\begin{equation} \label{eq: multivar mf dynamics}
    \dd X_t = \alpha_t \, \dd t + \sigma \, \dd W_t. 
\end{equation}
Recall that $m = \int_{\R^d} x \, \mu(\dd x)$.

\subsection{Mean Field Game}
Traditional methods for deriving the LQ problem solution begin with recovering the value function
\begin{equation}
\begin{split}
    v(x) \coloneqq \inf _{\alpha \in \mathbb{A}} \E \Biggl[ \int_0^\infty e^{-\beta t} \biggl(\frac{1}{2} \alpha_t^T \alpha_t + (X_t - C_2 m)^T C_1 (X_t - C_2 m)\\
    \left. \qquad {}+ (X_t - c_4)^T C_3 (X_t - c_4) + m^T C_5 m \biggr) \, \dd t \; \right\rvert X_0 = x \Biggr]
\end{split}
\end{equation}
as the solution of a Hamilton-Jacobi-Bellman (HJB) equation. In the MFG case, we denote the optimal value function as $\hat{v}$ and assume it takes the form of the ansatz
\begin{equation} \label{eq: multivar mfg ansatz}
    \hat{v}(x) = x^T \hat{\Gamma}_2 x + \hat{\Gamma}_1^T x + \hat{\Gamma}_0.
\end{equation}
The HJB equation for the value function in the asymptotic infinite horizon setting is
\begin{equation} \label{eq: multivar mfg hjb}
\begin{split}
    0 &= \beta \hat{v}(x) - H(x, \mu)\\
    & = \beta \hat{v}(x) - \inf_{\alpha} \left\{\alpha^T \nabla \hat{v}(x) + \frac{1}{2} \tr\left(\sigma \sigma^T \nabla^2 \hat{v}(x)\right) + \frac{1}{2} \alpha^T \alpha \right.\\
    &\quad + (x - C_2 m)^T C_1 (x - C_2 m) + (x - c_4)^T C_3 (x - c_4) + m^T C_5 m  \biggr\}
\end{split}
\end{equation}
where $H$ is the Hamiltonian for the given MFG problem and $\nabla^2$ is the Hessian operator. As the infimum term in \cref{eq: multivar mfg hjb} is quadratic in $\alpha$, the value at which the infimum is attained is
\begin{equation*}
\hat{\alpha}(x) = -\nabla \hat{v}(x) \quad \left(= -(\hat{\Gamma}_2 + \hat{\Gamma}_2^T) x - \hat{\Gamma}_1) \right).
\end{equation*}
Substituting this control into the HJB equation yields
\begin{align*}
    0 &= \beta \hat{v}(x) + \frac{1}{2} \nabla \hat{v}(x)^T \nabla \hat{v}(x) - \frac{1}{2} \tr\left(\sigma \sigma^T \nabla^2 \hat{v}(x)\right)\\
    &\quad - (x - C_2 m)^T C_1 (x - C_2 m) - (x - c_4)^T C_3 (x - c_4) - m^T C_5 m.
\end{align*}
Finally, we plug in our ansatz for the value function and simplify to obtain
\begin{align*}
    0 &= x^T (2 \hat{\Gamma}_2^2 + \beta \hat{\Gamma}_2 - C_1 - C_3) x\\
    & + \left((\beta I + 2 \hat{\Gamma}_2) \hat{\Gamma}_1 + 2 C_1 C_2 m + 2 C_3 c_4\right)^T x\\
    & + \beta \hat{\Gamma}_0 + \frac{1}{2} \hat{\Gamma}_1^T \hat{\Gamma}_1 - \tr(\sigma \sigma^T \hat{\Gamma}_2) - m^T (C_2^{T} C_1 C_2 + C_5) m - c_4^T C_3 c_4.
\end{align*}
From this we can solve for the coefficients in \cref{eq: multivar mfg ansatz} to get
\begin{align*}
    &\hat{\Gamma}_2 =\frac{1}{4}\left(-\beta I + [\beta^2 I + 8(C_1 + C_3)]^{1/2}\right)\\[.5em]
    &\hat{\Gamma}_1 = (\beta I + 2 \hat{\Gamma}_2)^{-1} \left(-2 C_1 C_2 m - 2 C_3 c_4 \right)\\[.5em]
    &\hat{\Gamma}_0 = \frac{1}{\beta}\left( -\frac{1}{2} \hat{\Gamma}_1^T \hat{\Gamma}_1 + \tr(\sigma \sigma^T \hat{\Gamma}_2) + m^T (C_2^T C_1 C_2 + C_5) m +c_4^T C_3 c_4 \right).
\end{align*}
Note that the square root of $\beta^2 I + 8(C_1 + C_3)$ in the first equation exists and is unique since this matrix is, in fact, symmetric positive-definite. As a result, $\hat{\Gamma}_2$ is also symmetric positive-definite which means we may rewrite the optimal control as
\begin{equation} \label{eq: multivar mfg control}
    \hat{\alpha}(x) = -2 \hat{\Gamma}_2 x - \hat{\Gamma}_1
\end{equation}

To recover the equilibrium mean field distribution, we will substitute \cref{eq: multivar mfg control} into the state SDE \cref{eq: multivar mf dynamics} which yields the Ornstein-Uhlenbeck process
\[
    \dd \hat{X}_t = -\left( 2\hat{\Gamma}_2 \hat{X}_t + \hat{\Gamma}_1 \right)\, \dd t + \sigma \, \dd W_t.
\]
Note that since $\hat{\Gamma}_2$ is symmetric positive-definite, $\hat{X}_t$ has a limiting distribution $\hat{\mu} = \lim_{t \to \infty} \mathcal{L}(\hat{X}_t)$ \cite{vatiwutipong2019}. Further, assuming $\hat{X}_0$ is normally distributed, the limiting distribution has an explicit form given by
\begin{equation}
    \hat{\mu} = \mathcal{N}\left( -\frac{1}{2}\hat{\Gamma}_2^{-1} \hat{\Gamma}_1, \hat{\Sigma} \right).
\end{equation}
Here, the covariance matrix $\hat{\Sigma}$ is defined by
\[
    \mathrm{vec}(\hat{\Sigma}) = \frac{1}{2}( \hat{\Gamma}_2 \oplus \hat{\Gamma}_2)^{-1} \mathrm{vec}(\sigma \sigma^T).
\]
where $\mathrm{vec}(A)$ is the vector obtained by stacking the columns of $A$ into a vector from left to right and $\oplus$ is the Kronecker sum \cite{vatiwutipong2019}.

Since the mean field interaction for the LQ problem is only through the mean $\hat{m} = \int x \,\hat{\mu}(\dd x)$, we note that a simplified form of $\hat{m}$ is
\begin{equation} \label{eq: multivar mfg mean}
    \hat{m} = -\frac{1}{2}\hat{\Gamma}_2^{-1} \hat{\Gamma}_1 = \left[C_1 + C_3 - C_1 C_2 \right]^{-1} C_3 c_4.
\end{equation}

\subsection{Mean Field Control}
As done previously, we denote the MFC value function by $v^*$ and claim that is has the form 
\begin{equation} \label{eq: multivar mfc ansatz}
v^*(x) = x^T \Gamma^*_2 x + {\Gamma_1^*}^T x + \Gamma^*_0.
\end{equation}
The MFC HJB equation differs from the MFG HJB equation (\ref{eq: multivar mfg hjb}) with the addition of an extra term involving the derivative of the Hamiltonian with respect to the measure $\mu$:
\begin{equation} \label{eq: multivar mfc hjb}
    0 = \beta v^*(x) - H(x, \mu) - \int_{\R^n} \frac{\delta H}{\delta \mu}(h, \mu) (x) \, \mu(\dd h).
\end{equation}
The derivative of $H$ with respect to $\mu$ in the sense of $L$-derivatives \cite{carmona-mfg} is
\begin{align*}
    \frac{\delta H}{\delta \mu}(h, \mu) (x) &= \frac{\delta}{\delta \mu} \left[ (h - C_2 m)^T C_1 (h - C_2 m) + m^T C_5 m\right] (x)\\[.5em]
    &= \frac{\delta}{\delta \mu} \left[ \left (h - C_2 \int_{\R^n} y \, \mu(\dd y) \right)^T C_1 \left (h - C_2 \int_{\R^n} y \, \mu(\dd y) \right) \right.\\
    & \quad \left. + \left(\int_{\R^n} y \, \mu(\dd y) \right)^T C_5 \int_{\R^n} y \, \mu(\dd y) \right] (x)\\[.5em]
    &= -2(h - C_2 m)^T C_1 C_2 x + 2 m^T C_5 x,
\end{align*}
and the integral in \cref{eq: multivar mfc hjb} is therefore
\[
    \int_{\R^n} \frac{\delta H}{\delta \mu}(h, \mu) (x) \, \mu(\dd h)
    = -2m^T (I-C_2)^T C_1 C_2 x + 2 m^T C_5 x.
\]

Recalling that the optimal control is given by
\begin{equation*} 
    \alpha^*(x) = -\nabla v^*(x) \qquad \biggl(= -(\Gamma^*_2 + {\Gamma^*_2}^T) x - \Gamma^*_1 \biggr),
\end{equation*}
we plug this into \cref{eq: multivar mfc hjb} and simplify to obtain
\begin{align*}
    0 &= x^T (2 {\Gamma^*_2}^2 + \beta \Gamma^*_2 - C_1 - C_3) x\\
    & + \left((\beta I + 2 \Gamma^*_2) \Gamma^*_1 + 2( C_1 C_2 + C_2^T C_1(I - C_2) - C_5) m + 2 C_3 c_4 \right)^T x\\
    & + \beta \Gamma^*_0 + \frac{1}{2} {\Gamma^*_1}^T \Gamma^*_1 - \tr(\sigma \sigma^T \Gamma^*_2) - m^T (C_2^{T} C_1 C_2 + C_5) m - c_4^T C_3 c_4.
\end{align*}
We solve for the coefficients of $v^*(x)$ as before to get
\begin{align*}
    \Gamma^*_2 &= \frac{1}{4}\left(-\beta I + [\beta^2 I + 8(C_1 + C_3)]^{1/2}\right)\\
    \Gamma^*_1 &= -2{(\beta I + 2 \Gamma^*_2)}^{-1} \left[ (C_1 C_2 + C_2^T C_1(I - C_2) - C_5) m + C_3 c_4 \right]\\
    \Gamma^*_0 &= \frac{1}{\beta}\left( -\frac{1}{2} {\Gamma^*_1}^T \Gamma^*_1 + \tr(\sigma \sigma^T \Gamma^*_2) + m^T (C_2^T C_1 C_2 + C_5) m +c_4^T C_3 c_4 \right).
\end{align*}
As before, we observe that $\Gamma^*_2$ must be symmetric positive-definite, so the optimal control can be rewritten as
\begin{equation} \label{eq: multivar mfc control}
    \alpha^*(x) = -2 \Gamma^*_2 x - \Gamma^*_1
\end{equation}

To derive the optimal mean field distribution, we again substitute \cref{eq: multivar mfc control} into \cref{eq: multivar mf dynamics} to get the Ornstein-Uhlenbeck process
\[
    \dd X^*_t = -\left(2 \Gamma_2^* X^*_t + \Gamma^*_1\right)\, \dd t + \sigma \, \dd W_t,
\]
whose limiting distribution $\mu^* = \lim_{t \to \infty} \mathcal{L}(X^*_t)$ is
\begin{equation}
    \mu^* = \mathcal{N}\left( -\frac{1}{2} {\Gamma^*_2}^{-1} \Gamma^*_1, \Sigma^* \right),
\end{equation}
assuming that $X^*_0$ is normally distributed.
The covariance matrix $\Sigma^*$ is again defined by
\[
    \mathrm{vec}(\Sigma^*) = \frac{1}{2}{( \Gamma^*_2 \oplus \Gamma^*_2)}^{-1} \mathrm{vec}(\sigma \sigma^T).
\]
Since the mean field interaction is only through the mean $m^* = \int x \, \mu^*(\dd x)$, we note that an equation for $m^*$ which only depends explicitly on the running cost coefficients is
\begin{equation}\label{eq: multivar mfc mean}
    m^* = \left[C_1 + C_3 - C_1 C_2 - C_2^T C_1 (I - C_2) + C_5 \right]^{-1} C_3 c_4.
\end{equation}

\end{document}